\input{epsf}

\documentclass[11pt,letterpaper]{article}

\usepackage[hang,perpage]{footmisc}

\usepackage{amsmath}
\usepackage{amsthm}
\usepackage{amssymb}
\usepackage{stmaryrd}
\usepackage{graphicx}
\usepackage{latexsym}
\usepackage{times}
\usepackage[varumlaut]{yfonts}
\usepackage[usenames]{color}
\usepackage{braket}
\usepackage{mathabx}
\usepackage{textcomp}
\usepackage{makeidx}
\usepackage{tocbibind}
\usepackage[bookmarks=true,backref=true,colorlinks=true,linkcolor=darkolivegreen,
citecolor=darkolivegreen,urlcolor=darkolivegreen]{hyperref}
\usepackage[T1]{fontenc}
\usepackage{mathrsfs}
\numberwithin{equation}{subsection}
\allowdisplaybreaks

\theoremstyle{definition}
\newtheorem{ass}{General Assumption}[subsection]
\newtheorem{theorem}[ass]{Theorem}
\newtheorem{lemma}[ass]{Lemma}
\newtheorem{definition}[ass]{Definition}
\newtheorem{corollary}[ass]{Corollary}

\newtheorem{ex}[ass]{Example}

\newtheorem{rem}[ass]{Remark}

\makeindex

\newcommand{\captionfonts}{\footnotesize}
\makeatletter  % Allow the use of @ in command names
\long\def\@makecaption#1#2{%
  \vskip\abovecaptionskip
  \sbox\@tempboxa{{\captionfonts #1: #2}}%
  \ifdim \wd\@tempboxa >\hsize
    {\captionfonts #1: #2\par}
  \else
    \hbox to\hsize{\hfil\box\@tempboxa\hfil}%
  \fi
  \vskip\belowcaptionskip}
\makeatother   % Cancel the effect of \makeatletter
\definecolor{darkolivegreen}{rgb}{0.333333, 0.419608, 0.1843140}

\makeatother

\begin{document}

% Uncomment for double-spacing \doublespacing

\title{Results on the diffusion equation with rough coefficients}

\author{Burak Aksoylu \\
Louisiana State University \\ 
Department of Mathematics \& \\
Center for Computation and Technology\\
Baton Rouge, LA 70803, USA
\and
Horst R. Beyer \\
Louisiana State University \\
Center for Computation and Technology \\
Baton Rouge, LA 70803, USA}

\date{\today}                                     
% Uncomment if TWO COLUMNS REQUIRED
%\twocolumn

\maketitle

\begin{abstract}

We study the behaviour of the solutions of the stationary diffusion
equation as a function of a possibly rough ($L^{\infty}$-)
diffusivity. This includes the boundary behaviour of the solution
maps, associating to each diffusivity the solution corresponding to
some fixed source function, when the diffusivity approaches infinite
values in parts of the medium. In $n$-dimensions, $n \geq 1$, by
assuming a weak notion of convergence on the set of diffusivities, we
prove the strong sequential continuity of the solution maps. In
$1$-dimension, we prove a stronger result, i.e., the unique
extendability of the map of solution operators, associating to each
diffusivity the corresponding solution operator, to a sequentially
continuous map in the operator norm on a set containing
`diffusivities' assuming infinite values in parts of the medium. In
this case, we also give explicit estimates on the convergence
behaviour of the map.
\end{abstract}

%Mathematics Subject Classification (2000)
%35J25 Boundary value problems for second-order, elliptic equations
%47F05 Partial differential operators
%65J10 Equations with linear operators
%65N99 None of the above, but in this section

\section{Motivation}

Numerical methods for the diffusion equation with rough coefficients
have been studied
extensively~\cite{ABDP:1981,bakhvalov1,BMMS_fosls:2005,KlWiDr:2002,knyazev1,
KnWi:2003,Oswald.P1999c} in the preconditioning (multigrid, domain
decomposition, and related iterative methods) literature starting the
early eighties and still continue to be an active area of research in
various preconditioning efforts~\cite{XuZh:2008,zhu:2008}.  This
article came about out of a need of deeper understanding of the
performance of preconditioners and their connection to the underlying
PDE.
\newline
\linebreak
In a recent article~\cite{AGKS:2007}, the first author constructed a
new preconditioning strategy with rigorous justification which is
comparable to algebraic multigrid.  It is shown in~\cite{AGKS:2007}
that analytical tools such as singular perturbation analysis gives
valuable insight about the asymptotic behavior of the solution of the
underlying PDE, hence, provides feedback for preconditioner
construction.
\newline
\linebreak
According to experience, the performance of a preconditioner depends
essentially on the degree to which the preconditioner operator
approximates the underlying operator.  Then, the fundamental need is
to explain the effectiveness of the preconditioner and to justify that
rigorously.  In that respect, one can view the tools in this article
as steps towards adding tools to the arsenal of methods of analysis
for rigorous justification at the interface of preconditioning
and operator theory.  Direct connections from the results here to
preconditioning will be the subject for future research.

\section{Introduction}
The diffusion equation
\begin{equation} \label{diffusioneq}
\frac{\partial u}{\partial t} = \textrm{div}  
\left(\, p \, \textrm{grad} \, u \right) + f  
\end{equation}
describes general diffusion processes, including the propagation 
of heat, and flows through porous media. Here $u$ is 
the density of the diffusing material, $p$ is the diffusivity 
of the material, and 
the function $f$ describes the distribution of `sources' and `sinks'.
This paper focuses on stationary solutions of 
(\ref{diffusioneq}) satisfying
\begin{equation} \label{diffusioneqstat}
- \textrm{div}  
\left(\, p \, \textrm{grad} \, u \right) = f  \, \, .
\end{equation}
For instance, the fictitious domain method and composite materials are
sources of rough coefficients; see the references in~\cite{KnWi:2003}.
Important current applications deal with composite materials whose
components have nearly constant diffusivity, but vary by several
orders of magnitude.  In composite material applications, it is quite
common to idealize the diffusivity by a piecewise constant function
and also to consider limits where the values of that function approach
zero or infinity in parts of the material.
\newline
\linebreak
Results of such study were given first by J. L. Lions~\cite{lions}.
In his lecture notes, he considers the limit of the solution of
(\ref{diffusioneqstat}) where the limit is associated to a
one-parameter family of piecewise constant diffusivities approaching
zero on a subdomain.  The same piecewise constant one-parametric
approach was used in~\cite{bakhvalov1,knyazev1}, but with
diffusivities approaching infinity on a subdomain. The limitation of
the one-parametric approach is its dependence on the particular
approximating sequence.  To the knowledge of the authors, this paper
is the first to address these questions in the necessary generality.
Hence, we consider general families of diffusivities that
are not necessarily piecewise constant.  In addition, due to the
atomistic structure of matter, the physical treatment of diffusion
involves regular ($C^\infty$-) diffusivity. It is unclear to what
extent the idealization of diffusivity by piecewise constant
coefficients has the capability to capture the underlying physics.
Mathematically, the severe contrast in diffusivity should be
represented by a regular function whose size is changing drastically
over small distances in interface regions. In this paper, we
demonstrate that the assumption of piecewise constant diffusivities is
meaningful by showing a continuous dependence of the solutions on the
diffusivity.
\newline
\linebreak
Furthermore, the diffusion equation is meaningless if the
`diffusivity' is infinite or zero in regions of the material.  Physics
requires nowhere vanishing diffusivity in the interior of the
material. As a consequence, only the relative size of diffusivities
should be significant.  Therefore, physically, one might expect that
both types of the above limits are equivalent, but mathematically
there are differences.  The limit of the solution as diffusivity
approaches infinite values exists.  However, only the limit of the
scaled solution exists as diffusivity approaches zero values (see
Example~\ref{counterex} for both cases). That is why, we choose to
work with diffusivity approaching infinity. We will refer these cases
as `asymptotic cases'.
\newline
\linebreak
Also, the treatment in \cite{bakhvalov1,knyazev1} considers only
limits on specific parts of the material.  In this connection, it
should also be remarked that, although (\ref{diffusioneq}),
(\ref{diffusioneqstat}) are linear equations, in general, their
solutions depend non-linearly on the coefficients.
\newline
\linebreak
For the treatment of these questions, we use methods from 
operator theory. For this, we use a common approach to 
give (\ref{diffusioneq}) a well-defined meaning that, in a 
first step, represents the diffusion 
operator 
\begin{equation} \label{diffusionoperator}
- \textrm{div} \, p \, \textrm{grad}  
\end{equation}
as a densely-defined positive self-adjoint linear operator 
$A_{p}$ in a suitable Hilbert space. As a result, (\ref{diffusioneqstat}) 
is represented by the equation
\begin{equation*}
A_{p} u = f \, \, , 
\end{equation*}
where $f$ is an element of the Hilbert space, and $u$ is 
from the domain, $D(A_{p})$, of $A_{p}$.
\footnote{After that, the abstract 
theory of strongly continuous one-parameter semigroups 
of operators can be used to associate a rigorous  
formulation of a well-posed  
initial value problem to (\ref{diffusioneq}) 
\cite{beyer,engel,pazy}. In this, $A_{p}$ becomes the infinitesimal
generator of time evolution. 
This last step will not be detailed
here.}  
\newline
\linebreak
Specifically, we treat the class ${\cal L}$ of 
diffusivities $p \in L^{\infty}(\Omega)$ that  
are almost everywhere $\geq \varepsilon$ 
on $\Omega$ for 
some $\varepsilon >0$, where $\Omega \subset {\mathbb{R}}^n$, 
$n \in {\mathbb{N}}^{*}$, is some non-empty 
open subset. By use of Dirichlet boundary conditions, 
it defines $A_{p}$ as an operator in the complex Hilbert space 
$L^2_{\mathbb{C}}(\Omega)$.
For non-smooth $p$, the domain of $A_{p}$ depends 
heavily on $p$. This fact significantly complicates the study of
sequences of functions of $A_{p}$.
\newline
\linebreak 
In this paper, we turn to a first-order formulation of
(\ref{diffusioneqstat}) which is often referred as mixed formulation
in the discretization literature~\cite{braess}. The first-order
formulation was popularized in the least squares finite element
community by the so-called FOSLS pioneering
paper~\cite{CLMM:1994_foslsFirstPaper}.  Here, we provide the
self-adjointness of the corresponding operator ${\hat{A}}_{p}$ in a
Hilbert space. The key property of ${\hat{A}}_{p}$ is that its domain,
$D({\hat{A}}_{p})$, is independent of $p$. This property is exploited
in establishing the continuity of the solutions $A_{p}^{-1} f$ as a
function of $p$.  Moreover, ${\hat{A}}_{p}$ remains defined for the
asymptotic cases when (\ref{diffusioneq}), (\ref{diffusioneqstat}) are
ill-defined. This fact is used in the study of the asymptotic cases.
\newline
\linebreak
Specifically, for $p \in {\cal L}$ and by assuming a weak notion of
convergence in ${\cal L}$, we show that the maps that associate $p$ to
the operator $A_{p}^{-1}$ and $- p \, \nabla A_{p}^{-1}$,
respectively, are strongly sequentially continuous, see
Theorem~\ref{stongconvergence} and
Corollary~\ref{stongconvergencecor}.  In particular, this shows in
these cases that the approximation by discontinuous coefficients to
physical diffusivity is indeed meaningful.  In addition, for the case
$n=1$ and bounded open intervals of ${\mathbb{R}}$, we show
stronger results that include also the asymptotic cases, except that
where the asymptotic `diffusivity' is almost everywhere infinite on
the interval.  In this case, the maps that associate ${\bar{p}}$ to
the operator $A_{1/{\bar{p}}}^{-1}$ and $- (1/{\bar{p}}) \, \nabla
A_{1/{\bar{p}}}^{-1}\,$, respectively, have unique extensions to
sequentially continuous maps in the operator norm on the set of
a.e. positive elements of $L^{\infty}(\Omega) \setminus \{0\}$, see
Corollary~\ref{uniformconvergence1dim},~\ref{uniformconvergence1dim2}.
In addition, an explicit estimate of the convergence behaviour of the
maps is given, see Theorem~\ref{speedofconv}. It is still an open
problem, whether the last results are generalizable to dimensions
$n \geq 2$.

\section{Basic notation}

\label{prerequisites}

Mainly, this section introduces basic notation. In particular, an
operator theoretic definition of Sobolev spaces is given that is based
on weak derivative operators, instead of distributions. In such
formulation, the completeness of the Sobolev spaces is an obvious
consequence of the cussedness of these operators.  Also, we give some
basic results that are connected to this formulation. For the
convenience to the reader, corresponding proofs are given in the
appendix.

\begin{ass}

In the following, let $n \in {\mathbb{N}}^{*}$ and $\Omega$ be a
non-empty open subset of ${\mathbb{R}}^n$.

\end{ass}

We follow common usage and do not differentiate between a function $f$
which is almost everywhere defined (with respect to a chosen measure)
on some set and the associated equivalence class consisting of all
functions which are almost everywhere defined on that set and differ
from $f$ only on a set of measure zero. The following definitions need
to be understood in this sense.

\begin{definition} ({\bf Complex $L^p$-spaces})
\begin{itemize}
\item[(i)]
For $p > 0$, the symbol
$L_{\mathbb{C}}^p(\Omega)$     
denotes the vector space of all complex-valued measurable functions
$f$ which are a.e. defined on $\Omega$ and such that 
$|f|^p$ is integrable  
with respect to the 
Lebesgue measure $v^n$. For every such $f$, 
we define the $L^p$-norm $\|f\|_p$ corresponding to $f$ by 
\begin{equation*}
\|f\|_p := \left( \, \int_{\Omega} |f|^p \, dv^n \right)^{1/p} \, \, .
\end{equation*}   
In addition, for the special case
$p=2$, 
we define a scalar product 
$\braket{\,|\,}_2$ 
on $L_{\mathbb{C}}^2(\Omega)$ by
\begin{equation*} 
\braket{f|g}_2 := \int_{\Omega} f^{*} g \, dv^n  \, \, ,
\end{equation*}
for all $f,g \in L_{\mathbb{C}}^2(\Omega)$. Here
$^*$ denotes complex 
conjugation on ${\mathbb{C}}$. As a consequence,
$\braket{\,|\,}_2$ is antilinear 
in the first argument and linear in its second. This convention 
will be used for sesquilinear forms in general.  
\item[(ii)]
$L_{\mathbb{C}}^{\infty}(\Omega)$ denotes the 
vector space of complex-valued measurable bounded functions
on $\Omega$. For every $f \in L_{\mathbb{C}}^{\infty}(\Omega)$, we 
define 
\begin{equation*}
\|f\|_{\infty} := \sup_{x \in \Omega}{|f(x)|} \, \, . 
\end{equation*}
\item[(iii)] For every $k \in {\mathbb{N}}^{*}$ and 
$f, g \in (L_{\mathbb{C}}^2(\Omega))^k$, we define 
\begin{equation*} 
\braket{f|g}_{2,k} := \sum_{j=1}^{k} \braket{f_{j}|g_{j}}_2  
\, \, , \, \, 
\|f\|_{2,k} := \left(\,\sum_{j=1}^{k} \|f_{j}\|_{2}^2 \right)^{1/2} 
\, \, . 
\end{equation*}
\end{itemize}
\end{definition}

\begin{definition} ({\bf Weak derivatives and Sobolev spaces})
We define
\begin{itemize}
\item[(i)]
for every multi-index $\alpha \in {\mathbb{N}}^n$
the densely-defined linear operator 
$\partial^{\, \alpha}$ in 
$L^{2}_{\mathbb{C}}(\Omega)$ 
\glossary{$\partial^{\, \alpha}$,\, $\alpha$-th weak 
partial derivative} 
by
\begin{equation*}
\partial^{\, \alpha}  \, := \,
(-1)^{|\alpha|}.    
\left( C^{\, \infty}_{0}(\Omega,{\mathbb{C}}) 
\rightarrow L^{2}_{\mathbb{C}}(\Omega), 
f \mapsto  \frac{\partial^{\, \alpha} f}{\partial x^\alpha} \right)^{*}
\, \, , 
\end{equation*}
where $*$ denotes the adjoint operation and  
\begin{equation*}
|\alpha| := \sum_{j=1}^{n} \alpha_{j} \, \, .
\end{equation*}

\item[(iii)]

for every $k \in {\mathbb{N}}$ 
\glossary{$W^{n}_{\mathbb{C}}(\Omega)$,\, Sobolev space of order $n$}
the Sobolev space $W^{k}_{\mathbb{C}}(\Omega)$ of order $k$ by

\begin{equation*}
W^{k}_{\mathbb{C}}(\Omega) := \bigcap_{\alpha \in {\mathbb{N}}^n, 
|\alpha| \leq k} D(\partial^{\, \alpha}) \, \, .
\end{equation*}  
Equipped with the scalar product 
\begin{equation*}
\braket{\, , \,}_{k} \, \, \, : 
\, (W^{k}_{\mathbb{C}}(\Omega))^2 \rightarrow {\mathbb{C}}
\, \, , 
\end{equation*}
defined by \glossary{$\braket{\, , \,}_{k}$}
\begin{equation*}
\braket{f,g}_{k} \, := \, \sum_{\alpha \in {\mathbb{N}}^n,|\alpha| \leq \, k} 
\braket{\partial^{\, \alpha}f|\partial^{\, \alpha}g}_{2} 
\end{equation*}
for all $f,g \in W^{k}_{\mathbb{C}}(\Omega)$, 
$W^{k}_{\mathbb{C}}(\Omega)$ becomes a Hilbert space.

\item[(iv)]
$W^k_{0,{\mathbb{C}}}(\Omega)$ as the closure of 
$C_{0}^{\infty}(\Omega,{\mathbb{C}})$ in 
$(W^k_{{\mathbb{C}}}(\Omega),\vvvert \, \vvvert_{k})$, where 
$\vvvert \, \vvvert_{k}$ denotes the norm that is induced 
on $W^k_{{\mathbb{C}}}(\Omega)$ 
by $\braket{\,,\,}_{k}$.

\end{itemize}

\end{definition}

We note that 

\begin{lemma} \label{partialintegration} ({\bf Partial integration})
\begin{equation} \label{partialintegration3}
\braket{f|\partial^{\,e_k} g}_2 = - \braket{\partial^{\,e_k} f|g}_2  
\end{equation}
for all $(f,g) \in  
W^1_{0,{\mathbb{C}}}(\Omega) \times  W^1_{{\mathbb{C}}}(\Omega)$
and $k \in {\mathbb{N}}^{*}$, where $e_k$ denotes the $k$-th 
canonical unit vector of ${\mathbb{R}}^{n}$.
\end{lemma}

The next defines gradient operators. 

\begin{definition} 
({\bf Gradient operators})
We define the $(L^{2}_{\mathbb{C}}(\Omega))^{n}$-valued
densely-defined linear operators in $L^{2}_{\mathbb{C}}(\Omega)$
\begin{align*}
& \nabla_{\!0} : C^{\infty}_{0}(\Omega,{\mathbb{C}}) \rightarrow 
(L^{2}_{\mathbb{C}}(\Omega))^{n} \, \, , \, \, 
\nabla_{\!w} : W^1_{{\mathbb{C}}}(\Omega) \rightarrow 
(L^{2}_{\mathbb{C}}(\Omega))^{n}
\end{align*}
by  
\begin{align*}
& \nabla_{0} f := 
{\phantom{\bigg|}}^{\!t}\!\left(\frac{\partial f}{\partial x_1},  
\dots, \frac{\partial f}{\partial x_n}\right) \, \, , \, \, 
\nabla_{\!w} g := {\phantom{}}^{t}(\partial^{\,e_1}g,\dots, 
\partial^{\,e_n}g)
\end{align*}
for all $f \in C^{\infty}_{0}(\Omega,{\mathbb{C}})$ and 
$g \in W^1_{{\mathbb{C}}}(\Omega)$.
\end{definition}

Then the following holds.

\begin{lemma} \label{gradientoperatorsandadjoints}
({\bf Adjoints of gradient operators})
\begin{equation} \label{gradientsandadjoints}
({\nabla_{\!0}}^{*})^{*} = \nabla_{\!w}
\big|_{W^{1}_{0, \mathbb{C}}(\Omega)}
\, \, , \, \,
\left( \nabla_{\!w}
\big|_{W^{1}_{0, \mathbb{C}}(\Omega)}\right)^{*} = 
{\nabla_{\!0}}^{*} \, \, .
\end{equation}
\end{lemma}

\section{Basic properties of the diffusion operator}

This section provides the basis of the paper. It defines 
the diffusion operator 
as operator in $L^{2}_{\mathbb{C}}({\Omega})$ and gives 
basic properties.  

\begin{definition} \label{defA} Let $\bar{p} : 
\Omega \rightarrow {\mathbb{R}}$ be measurable and such that 
$1 / \bar{p}$ is a.e. defined on $\Omega$. We define the linear 
operator $A : D(A) \rightarrow L^{2}_{\mathbb{C}}({\Omega})$ 
in $L^{2}_{\mathbb{C}}({\Omega})$ by 
\begin{equation*}
D(A) := \{ u \in W^1_{0,{\mathbb{C}}}({\Omega}): 
(1 / \bar{p}) {\nabla_{w}} u \in D({\nabla_{\!0}}^{*}) \}
\end{equation*}
and 
\begin{equation*}
A u := {\nabla_{\!0}}^{*} (1 / \bar{p}) \, {\nabla_{w}} u
\end{equation*} 
for every $u \in D(A).$
\end{definition}

Diffusion operators corresponding to diffusivities from 
the following large subset ${\cal L}$ of $L^{\infty}(\Omega)$ 
will turn out to be densely-defined self-adjoint linear operators. 

\begin{definition}
We define the subset ${\cal L}$ 
of $L^{\infty}(\Omega)$ to consist 
of those elements $\bar{p}$ for which there are real 
$C_1, C_2$ satisfying $C_2 \geq C_1 > 0$ and such that
$C_1 \leq \bar{p} \leq C_2$ a.e. on $\Omega$. Note that the last 
also implies that $1 / \bar{p} \in {\cal L}$ 
and in particular that  
$1 / C_2 \leq 1 / \bar{p} \leq 1/ C_1$ a.e. on $\Omega$. 
\end{definition}

The next proves the self-adjointness of diffusion operators 
corresponding to diffusivities from ${\cal L}$. For this,
so called `form methods' from operator theory are used. For 
these methods, see 
\cite{kato}.  

\begin{theorem} \label{secondorderoperator}
Let $\bar{p} \in {\cal L}$.
Then $A$ is a densely-defined 
linear self-adjoint operator in $L^{2}_{\mathbb{C}}({\Omega})$. 
\end{theorem}

\begin{proof}
For this, we define a positive Hermitian sesquilinear form 
$s : (W^1_{0,{\mathbb{C}}}({\Omega}))^2 \rightarrow
{\mathbb{C}}$ by 
\begin{equation*}
s(u,v) := \braket{{\nabla_{w}} u \, | \, (1/\bar{p}) \,  
{\nabla_{w}} v}_{2,n}  
\end{equation*}
for all $u,v \in W^1_{0,{\mathbb{C}}}({\Omega})$. Then 
$\braket{\,|\,}_{s} : (W^1_{0,{\mathbb{C}}}({\Omega}))^2 \rightarrow
{\mathbb{C}}$, defined by 
\begin{equation*}
\braket{u|v}_{s} := s(u,v) + \braket{u|v}_{2}
\end{equation*}
for every $u,v \in W^1_{0,{\mathbb{C}}}({\Omega})$, defines a 
scalar product on $W^1_{0,{\mathbb{C}}}({\Omega})$ with induced 
norm $\|\,\,\|_{s} : W^1_{0,{\mathbb{C}}}({\Omega}) \rightarrow 
{\mathbb{R}}$ given by 
\begin{equation*}
\|u\|_{s}^2 =  \braket{{\nabla_{w}} u \, | \, (1/\bar{p}) \,  
{\nabla_{w}} u}_{2,n} + \|u\|_{2}^2 
\end{equation*} 
for all $u \in W^1_{0,{\mathbb{C}}}({\Omega})$. In particular, $s$
is closable. For the proof, let $u_1,u_2,\dots$ be a Cauchy sequence 
in $(W^1_{0,{\mathbb{C}}}({\Omega}),\|\, \,\|_{s})$ and such that 
\begin{equation*}
\lim_{\nu \rightarrow \infty} \|u_{\nu}\|_{2} = 0 \, \, .
\end{equation*}
We note that 
\begin{align*}
& \min\{1,1/C_2\} \, \vvvert u \vvvert_{1}^2 \leq 
\frac{1}{C_2} \, 
\|{\nabla_{w}} u \|_{2,n} + \|u\|_{2}^2  \leq 
\|u\|_{s}^2 \\
& \leq \frac{1}{C_1} \, 
\|{\nabla_{w}} u \|_{2,n} + \|u\|_{2}^2 
\leq \max\{1,1/C_1\} \, \vvvert u \vvvert_{1}^2 \, \, ,
\end{align*}
where $C_1,C_2 \in {\mathbb{R}}$ satisfy $C_2 \geq C_1 > 0$ 
and are such 
that $C_1 \leq \bar{p} \leq C_2$ a.e. on $\Omega$, 
and hence that $\|\,\,\|_{s}$ and the restriction of 
$\vvvert \, \, \vvvert_{_{1}}$ to $W^1_{0,{\mathbb{C}}}({\Omega})$
are equivalent. Hence it follows that
\begin{equation*}
\lim_{\nu \rightarrow \infty} \|u_{\nu}\|_{s} = 0 \, \, .
\end{equation*}
Since $(W^1_{0,{\mathbb{C}}}({\Omega}),\|\, \,\|_{s})$ is in particular
complete, it follows that $s$ coincides with its closure. As a 
consequence, there is a unique densely-defined linear self-adjoint 
operator $A : D(A) \rightarrow L^2_{{\mathbb{C}}}({\Omega})$ 
in $L^2_{{\mathbb{C}}}({\Omega})$ such that $D(A)$ is a dense subspace 
of $(W^1_{0,{\mathbb{C}}}({\Omega}), \vvvert \,\, \vvvert_{1})$
and such that 
\begin{equation*}
\braket{u|A u}_{2} = s(u,u) = 
\braket{{\nabla_{w}} u \, | \, (1/\bar{p}) \,  
{\nabla_{w}} u}_{2,n} 
\end{equation*}  
for all $u \in D(A)$. In particular, $D(A)$ consists 
of all $ u \in W^1_{0,{\mathbb{C}}}({\Omega})$ for which there
is $f \in  L^2_{{\mathbb{C}}}({\Omega})$ such that 
\begin{equation*}
\braket{\, f\,|\dots\,}_{2}\big|_{W^1_{0,{\mathbb{C}}}({\Omega})} = 
\braket{ \, (1/\bar{p}) \,  {\nabla_{w}} u \, | \, 
{\nabla_{w}} \dots \, }_{2,n}\!\big|_{W^1_{0,{\mathbb{C}}}({\Omega})} 
\, \, .
\end{equation*}  
Further, if $u$ and $f$ satisfy these requirements, then 
\begin{equation*}
Au = f \, \, .
\end{equation*}
Hence $u \in D(A)$ if and only if 
\begin{equation*}
(1/\bar{p}) \,  {\nabla_{w}} u \in D \left( \left( \nabla_{\!w}
\big|_{W^{1}_{0, \mathbb{C}}(\Omega)}\right)^{*} \, \right) = 
D({\nabla_{\!0}}^{*})
\end{equation*} 
and in this case 
\begin{equation*}
A u = {\nabla_{\!0}}^{*} \, (1/\bar{p}) \,  {\nabla_{w}} u 
\, \, . 
\end{equation*}
\end{proof}

For completeness, the next gives the proof that 
diffusion operators corresponding to diffusivities
from ${\cal L}$ have a purely discrete spectrum, i.e., 
that their spectrum is a discrete subset of the real numbers 
consisting of eigenvalues of finite multiplicity
and that there is a Hilbert basis consisting of eigenvectors. 
This result is not used in the following.

\begin{corollary} Let
$\bar{p} \in {\cal L}$ and, in addition, $\Omega$ be bounded.   
Then $A$ has a purely discrete spectrum.
\end{corollary}

\begin{proof}
According to the proof of Theorem~\ref{secondorderoperator},
$\|\, \,\|_{s} : W^1_{0,{\mathbb{C}}}({\Omega}) \rightarrow 
{\mathbb{R}}$ defines a norm which is equivalent to the restriction 
of $\vvvert \,\, \vvvert_{1}$ to $W^1_{0,{\mathbb{C}}}({\Omega})$.
Hence the closed unit ball $B$ in 
$(W^1_{0,{\mathbb{C}}}({\Omega}),\|\, \,\|_{s})$ is contained 
in a closed ball of 
$(W^1_{0,{\mathbb{C}}}({\Omega}),\vvvert \,\, \vvvert_{1})$. The last 
is relatively compact in $L^2_{{\mathbb{C}}}({\Omega})$. From this, it 
follows also that $B$ is relatively compact in 
$L^2_{{\mathbb{C}}}({\Omega})$. Hence it follows, see, e.g.,
\cite{reed} Vol. IV, that $A$ has a purely 
discrete spectrum.
\end{proof}

\begin{figure} 
\centering
\includegraphics[width=5.6cm,height=5.6cm]{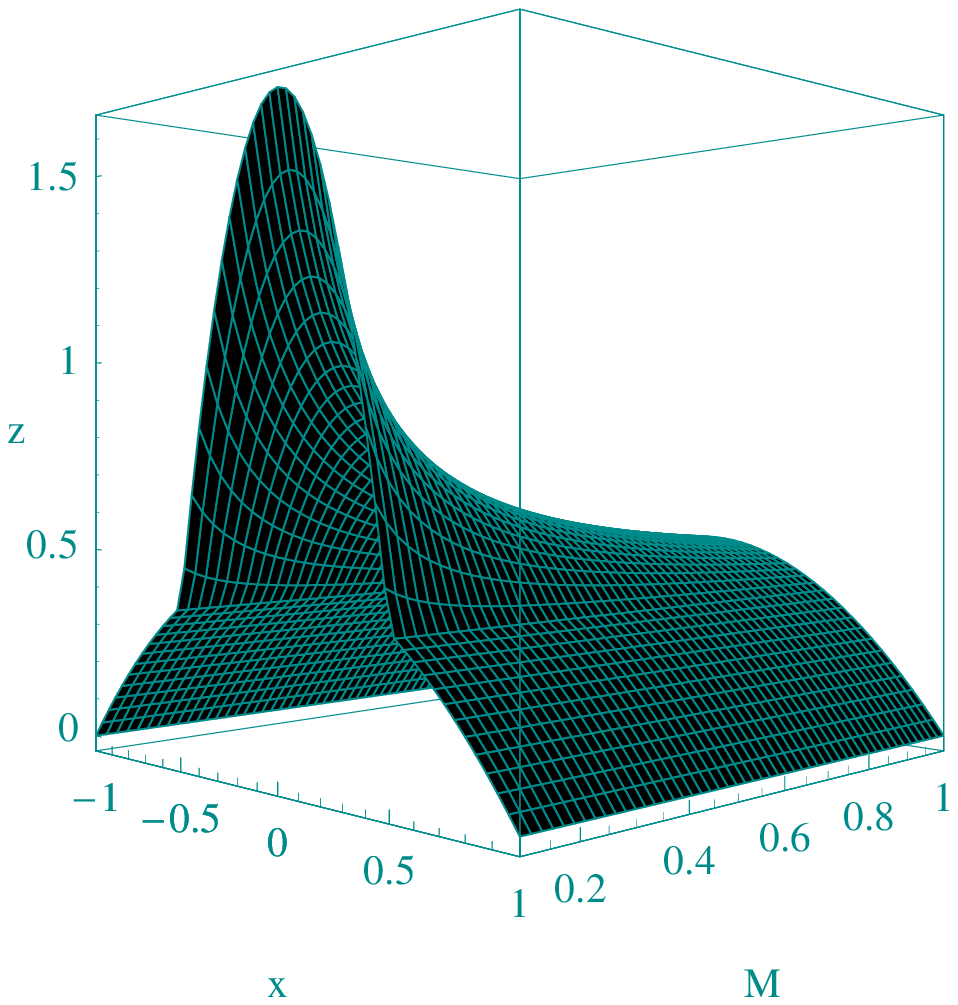} 
\quad 
\includegraphics[width=5.6cm,height=5.6cm]{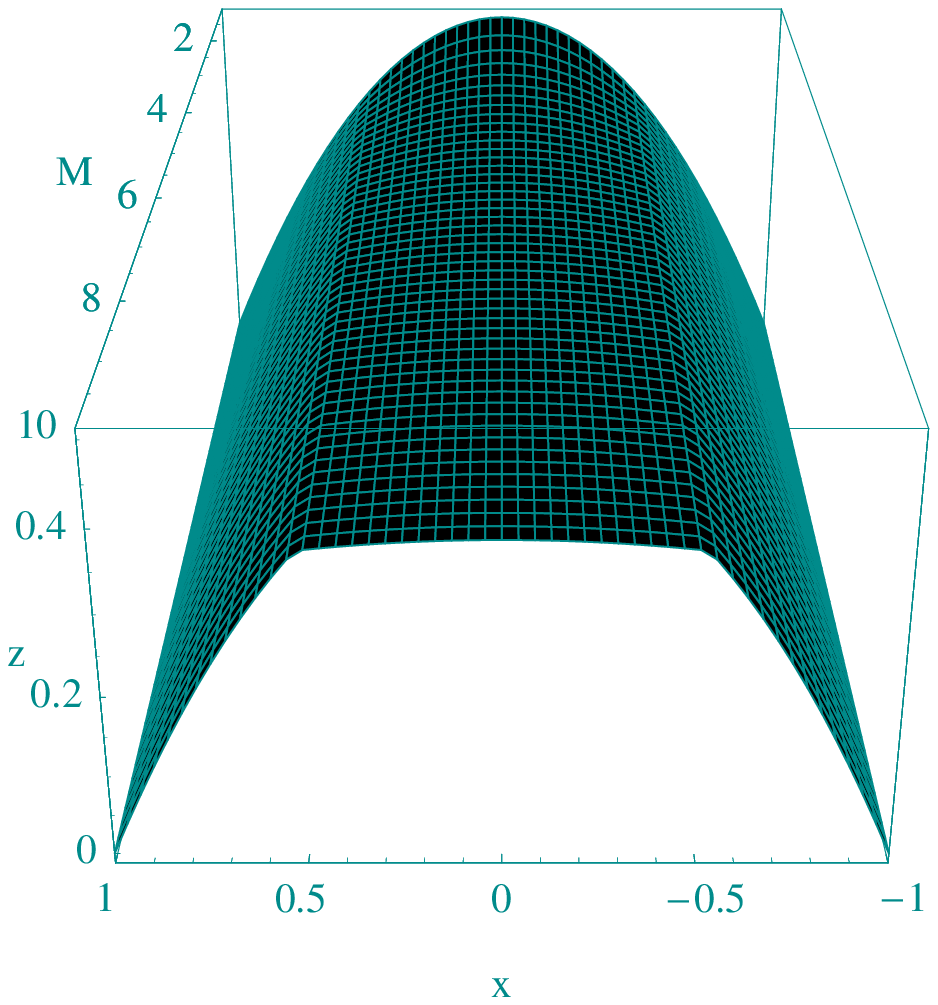}
\caption{Graphs of $u$ from Example~\ref{counterex} as a function 
of $M$.} 
\label{fig3}
\end{figure}

\begin{ex} \label{counterex} 
The following example illustrates the influence of discontinuities of 
the diffusivity on the regularity of the elements in $D(A)$. 
Consider the case that $\Omega = I :=
(-1,1)$ and a piecewise constant diffusivity $p : I \rightarrow
{\mathbb{R}}$ given by
\begin{equation*}
p(x) := \begin{cases}
\, \, \, 1 &\text{if $-1 < x < -1/2$}  \\
M &\text{if $-1/2 \leq x \leq 1/2$} \\
\, \, 1 & \text{if $1/2 < x < 1$}  
\end{cases}
\end{equation*}
for $x \in I$, where $M > 0$. Then $A u = f$, where 
$u : I \rightarrow {\mathbb{R}}$ is defined 
by 
\begin{align*}
u(x) := 
\begin{cases}
(1-x^2)/2 & \text{if $-1 < x \leq - 1/2$}  \\
(1 - 4 x^2 + 3 M )/(8 M)
& \text{if $-1/2 < x < 1/2$} \\
(1-x^2)/2 
& \text{if $1/2 \leq x < 1$}
\end{cases}
\end{align*}
and $f$ is the constant function on $I$ of value $1$. We note 
that $u^{\, \prime}$ has no extension to a continuous function
on $I$ if $M \neq 1$. In general, discontinuities in the diffusivity cause
low regularity of elements in $D(A)$. Also, see the concluding remarks.
\end{ex}

There is a unique solution $u_{f}$ to the equation
\begin{equation*}
A u_{f} = f
\end{equation*}
for every $f \in L^2_{{\mathbb{C}}}({\Omega})$ if and only 
if $A$ is bijective or equivalently, if and only if 
$0$ is not part of the spectrum of $A$. In general, $A$ is not 
bijective. For instance, the  
operator $A$ that is associated to $\Omega = {\mathbb{R}}^n$
and the diffusivity $p(x) = 1$ for every $x \in {\mathbb{R}}^n$
is not surjective. 
Below, we place a restriction on $\Omega$ that leads to bijective 
diffusion operators.

\begin{ass}  \label{poincarelemma}
In the following, we assume that 
$\Omega$ is in addition such that the following Poincare inequality 
is valid  
\begin{equation} \label{poincare}
 \|\, \partial^{e_{j}} f \, \|_2 \geq c \, \|f\|_2
\end{equation} 
for some $j \in \{1,\dots,n\}$ and every 
$f \in W^{1}_{0, \mathbb{C}}(\Omega)$, where $c > 0$. In the 
remainder, such $c$ is considered chosen.    
\end{ass}

\begin{rem}
It is known that $\Omega$ of the assumed type are not 
necessarily bounded. For instance, every non-trivial open set,
for which there is ${\bf n} \in {\mathbb{R}}^n \setminus \{ 0 \}$
along with real numbers 
$a, b$ such that   
\begin{equation*}
a < x \cdot {\bf n} < b 
\end{equation*}
for all $x \in \Omega$, is of this type.  
\end{rem}

In particular, the following proves that diffusion operators 
corresponding to diffusivities from ${\cal L}$ are bijective. 

\begin{theorem} \label{Abijective} 
Let $\bar{p} \in {\cal L}$.
The spectrum $\sigma(A)$
of $A$ satisfies
\begin{equation} \label{specest} 
\sigma(A) \subset 
[\,c^2 / \, C ,\infty ) \, \, ,
\end{equation}
where $j \in \{1,\dots,n \}$ is such that 
${\bf n}_{j} \neq 0$ and $C > 0$ is such that 
$\bar{p} \leq C$ a.e. on $\Omega$.  
\end{theorem}

\begin{proof}
For this, let 
$j \in \{1,\dots,n \}$ be such that ${\bf n}_{j} \neq 0$.
For $u \in D(A)$, it follows that 
\begin{equation*}
\braket{u|A u}_{2} = 
\braket{{\nabla_{w}} u \, | \, (1/\bar{p}) \,  
{\nabla_{w}} u}_{2,n} \geq C^{-1}
\, \|{\nabla_{w}} u \, \|_{2,n}^2 \geq 
c^2 \, C^{-1} \, \| u \|_{2}^2 \, \, ,
\end{equation*}
where $C > 0$ is such that
$\bar{p} \leq C$ a.e. on $\Omega$.  
Hence it follows the validity of (\ref{specest}).
\end{proof}

\section{Properties of a first order operator connected to the diffusion
operator}

As indicated by Example~\ref{counterex}, for non-smooth diffusivities
$p$, the condition that $p \nabla_{w} u \in W^1_{\mathbb{C}}(\Omega)$
in the definition of the domain of $A$ leads to a strong dependence of
that domain on the diffusivity. This fact poses an obstacle to the
study of the map, associating to every diffusivity $p \in {\cal L}$
the corresponding operator $A^{-1}$, by the notion of strong resolvent
convergence, see,~\cite[Volume I, Section VIII.7]{reed},
~\cite[Section VIII, \S 1]{kato}. By use of the following vector
partial differential operator of the first order $\hat{A}$, this
problem can be circumvented. Its domain is independent of the
diffusivity. The connection of the resolvents of $A$ and $\hat{A}$ is
given in Theorem~\ref{Ahatbijective}.
  
\begin{definition} Let $\bar{p} \in L^{\infty}({\Omega})$. 
We define 
the densely-defined, linear 
operator
$\hat{A} :  W^1_{0,{\mathbb{C}}}({\Omega}) \times D({\nabla_{\!0}}^{*})
\rightarrow 
L^2_{\mathbb{C}}({\Omega}) \times 
(L^2_{\mathbb{C}}({\Omega}))^{n}$ 
in $L^2_{\mathbb{C}}({\Omega}) \times 
(L^2_{\mathbb{C}}({\Omega}))^{n}$ by 
\begin{equation*}
\hat{A} (u,q) := \left( \, 
{\nabla_{\!0}}^{*} q \, , \, {\nabla_{\!w}} u - {\bar{p}}
\, q \, \right) 
\end{equation*}
for every $(u,q) \in  W^1_{0,{\mathbb{C}}}({\Omega}) \times 
D({\nabla_{\!0}}^{*})$. 
\end{definition}

\begin{theorem}
The operator $\hat{A}$ is self-adjoint.
\end{theorem}

\begin{proof}
The statement is a consequence of 
Lemma~\ref{gradientoperatorsandadjoints}.
\end{proof}

The following gives a characterization of the kernel of 
$\hat{A}$. In particular, the result implies that 
$\hat{A}$ is bijective 
for diffusivities from ${\cal L}$.   

\begin{theorem} \label{characterizationkerAhat}
Let $\bar{p} \in L^{\infty}({\Omega})$ be a.e. positive. Then 
\begin{equation*}
\ker \hat{A}  =  \{0\} \times (\, \ker 
{\nabla_{\!0}}^{*} \cap \ker T_{\bar{p}}^n \, ) \, \, ,
\end{equation*}
where $T_{\bar{p}} \in L(L^{2}_{\mathbb{C}}({\Omega}),L_{\mathbb{C}}^{2}({\Omega}))$ denotes 
the maximal multiplication operator in 
$L^{2}_{\mathbb{C}}({\Omega})$ that 
is associated 
to $\bar{p}$.
\end{theorem}

\begin{proof}
'$\subset$': Let $q \in \ker {\nabla_{\!0}}^{*} \cap \ker T_{\bar{p}}^n$.
Then $(0,q) \in D(\hat{A})$ and 
\begin{equation*}
{\nabla_{\!0}}^{*} q  = 0 \, , \, - {\bar{p}}
\, q = 0 \, \, .
\end{equation*}  
Hence it follows that $(0,q) \in \ker \hat{A}$. \newline
'$\supset$':
Let $(u,q) \in \ker \hat{A}$. Then 
\begin{equation} \label{kernel}
{\nabla_{\!0}}^{*} q  = 0 \, , \, {\nabla_{\!w}} u - {\bar{p}}
\, q = 0 
\end{equation}
and hence 
\begin{align*}
& 0 = \braket{\,q \,| \, {\nabla_{\!w}} u - {\bar{p}} \, q \, }_{2,n} 
= \braket{\,q \,| \, {\nabla_{\!w}} u \, }_{2,n}
- \braket{\,q \,| \, {\bar{p}} \, q \, }_{2,n} \\
& = \braket{\,{\nabla_{0}}^{*} q \,| \, u \, }_{2} - 
\|\, {\bar{p}}^{1/2} q \, \|_{2,n} =  - 
\|\, {\bar{p}}^{1/2} q \, \|_{2,n}
\, \, .
\end{align*}
The last implies that 
\begin{equation*}
q \in \ker T_{{\bar{p}}^{1/2}}^n \, \, , \, \, 
\end{equation*}
where $T_{{\bar{p}}^{1/2}} 
\in L(L^{2}_{\mathbb{C}}({\Omega}),L^{2}_{\mathbb{C}}({\Omega}))$ denotes 
the maximal multiplication operator in $L^{2}_{\mathbb{C}}({\Omega})$ that 
is associated 
to ${\bar{p}}^{1/2}$, and hence also that 
\begin{equation*}
q \in \ker T_{{\bar{p}}}^n \, \,  . 
\end{equation*}
Further, by (\ref{kernel})-2), it follows that 
\begin{equation*}
{\nabla_{\!w}} u = 0 \, \, .
\end{equation*} 
The last implies that 
\begin{equation*}
{\nabla_{0}}^{*} {\nabla_{\!w}} u = 0 
\end{equation*}
and hence by Theorem~\ref{Abijective} that $u = 0$.
\end{proof}

The following example shows that the kernel of $\hat{A}$
is non-trivial if $\bar{p}$ vanishes on some open subset 
of $\Omega$. The vanishing of $\bar{p}$ on non-empty
subsets of 
$\Omega$ corresponds to the 
asymptotic cases mentioned in the introduction. 

\begin{ex} In the following, we give $q \in 
C_{0}^{\infty}({\mathbb{R}}^n,{\mathbb{R}}^n) 
\cap \ker {\nabla_{\!0}}^{*}$
for $n \geq 2$. For this, let $h$ be an element of  
$C^{\infty}_{0}({\mathbb{R}})$ with 
support contained in $[-1,1]$. In addition, let 
$\alpha$ be a non-zero antisymmetric $n \times n$-matrix. We 
define $q  \in  C^{\infty}_{0}({\mathbb{R}}^n,{\mathbb{R}}^n)$
by  
\begin{equation*}
q(x) := \frac{h(|x|^2)}{2} \sum_{i,j =1}^n \alpha_{ij} x_{j} 
e_{i}
\end{equation*}
for all $x = (x_1,\dots,x_n) \in {\mathbb{R}}^n$. Then  
\begin{equation*}
\sum_{i=1}^{n} \frac{\partial q_{i}}{\partial x_i}(x) =
h^{\, \prime}(|x|^2) \! \sum_{i,j=1}^{n} \alpha_{ij} x_{i} x_j 
= 0
\end{equation*} 
for all $x = (x_1,\dots,x_n) \in {\mathbb{R}}^n$ and hence 
$q \in \ker {\nabla_{\!0}}^{*}$.
\end{ex}

The following lemma prepares the subsequent theorem which estimates
the size of the gap around $0$ in the spectrum of $A$ and gives a
representation of the resolvent of $\hat{A}$ in terms of the resolvent
of $A$, i.e., (\ref{resolventCharacterization}). The main tool in
the proof is the Closed Graph Theorem in the form of see Theorem~3.1.9
in \cite{beyer}.

\begin{lemma} \label{auxiliaryoperators}
Let $\bar{p} \in {\cal L}$, $\sigma(A)$ the spectrum of $A$, 
$\lambda < \min \{\sigma(A)\}$ and $A_{\lambda} := A -\lambda$. 
Then 
\begin{itemize}
\item[(i)] ${\nabla_{\!w}} A_{\lambda}^{-1} \in L(L^{2}_{\mathbb{C}}({\Omega}),
(L^{2}_{\mathbb{C}}({\Omega}))^n)$\,,
\item[(ii)] $\overline{A_{\lambda}^{-1} {\nabla_{0}}^{*}} = 
(\,{\nabla_{\!w}} A_{\lambda}^{-1}\,)^{*}$,
\item[(iii)] $D(A_{\lambda}^{1/2}) = 
W^1_{0,{\mathbb{C}}}({\Omega})$ and $A_{\lambda}^{1/2} : 
W^1_{0,{\mathbb{C}}}({\Omega}) \rightarrow 
L^{2}_{\mathbb{C}}({\Omega})$ is continuous, 
\item[(iv)] $\overline{{\nabla_{\!w}} A_{\lambda}^{-1} {\nabla_{0}}^{*}
}$ is a positive self-adjoint element of 
$L((L^{2}_{\mathbb{C}}({\Omega}))^{n},(L^{2}_{\mathbb{C}}({\Omega}))^{n})$. 
\end{itemize}
\end{lemma}

\begin{proof}
`(i)': Since $\lambda \in {\mathbb{C}} \, \setminus \,  
\sigma(A)$, $A_{\lambda}$ is 
densely-defined, linear and bijective. 
Further, $A_{\lambda}$ is self-adjoint and 
strictly positive. As a consequence of its 
self-adjointness, $A_{\lambda}$ is in particular closed. Further, 
according to Lemma~\ref{gradientoperatorsandadjoints}, 
the restriction of ${\nabla_{\!w}}$ to $W^{1}_{0,{\mathbb{C}}}({\Omega})$
is a closed linear operator in $L^{2}_{\mathbb{C}}({\Omega})$ with values 
in $(L^{2}_{\mathbb{C}}({\Omega}))^n$. Since $D(A) \subset 
W^{1}_{0,{\mathbb{C}}}({\Omega})$
it follows from the closed graph theorem, e.g., see Theorem~3.16 
in \cite{beyer}, the existence of $C \in [0,\infty)$ such that 
\begin{equation*}
\|{\nabla_{\!w}} f \|_{2,n} \leq C \, \|A_{\lambda} f \|_{2}
\end{equation*}
for all $f \in D(A)$. As a consequence, it follows for every 
$f \in L^{2}_{\mathbb{C}}({\Omega})$ that 
\begin{equation*}
\|{\nabla_{\!w}} A_{\lambda}^{-1} f \|_{2,n} \leq C \, 
\|A_{\lambda} A_{\lambda}^{-1} f \|_{2}
= C \, \|f \|_{2} 
\end{equation*}
and hence that ${\nabla_{\!w}} A_{\lambda}^{-1} \in L(L^{2}_{\mathbb{C}}({\Omega}),
(L^{2}_{\mathbb{C}}({\Omega}))^n)$. \newline
`(ii)': $A_{\lambda}^{-1} {\nabla_{0}}^{*}$ is a densely-defined, linear
operator in $(L^{2}_{\mathbb{C}}({\Omega}))^n$ with values in $L^{2}_{\mathbb{C}}({\Omega})$.
Further, it follows for 
$f \in L^{2}_{\mathbb{C}}({\Omega})$ that 
\begin{equation*}
\braket{\, f \, | \, A_{\lambda}^{-1} {\nabla_{0}}^{*} q\,}_{2} = 
\braket{\,  A_{\lambda}^{-1} \! f \, | \, {\nabla_{0}}^{*} q\,}_{2} =
\braket{\,{\nabla_{\!w}} A_{\lambda}^{-1} f\, | \, q \, }_{2,n} 
\end{equation*}
for every $q \in D({\nabla_{0}}^{*})$ and hence that 
$f \in D((A_{\lambda}^{-1} {\nabla_{0}}^{*})^{*})$ as well as that 
\begin{equation*}
(A_{\lambda}^{-1} {\nabla_{0}}^{*})^{*} f = {\nabla_{\!w}} A_{\lambda}^{-1} f \, \, .
\end{equation*} 
As a consequence, 
\begin{equation*}
(A^{-1}_{\lambda} {\nabla_{0}}^{*})^{*} = {\nabla_{\!w}} 
A^{-1}_{\lambda} \, \, .
\end{equation*}
In particular, $A^{-1}_{\lambda} {\nabla_{0}}^{*}$ is closable and 
\begin{equation*}
\overline{A^{-1}_{\lambda} {\nabla_{0}}^{*}} = 
(\,{\nabla_{\!w}} A^{-1}_{\lambda}\,)^{*}  \, \, . 
\end{equation*}
`(iii)': In a first step, we prove the statement for the case $\lambda = 
0$. For this, we note that, as a consequence of 
Theorem~\ref{Abijective}, $0 < \min\{\sigma(A)\}$. Further, 
we note that $D(A)$ is a core $A^{1/2}$. For 
instance, this follows by Theorem~3.1.9 in \cite{beyer}. Hence 
$D(A)$ is dense in the Banach space 
$(D(A^{1/2}), \|\,\|_{A^{1/2}})$, where 
\begin{equation*}
\|f\|_{A^{1/2}} := \big[ \, \|f\|_{2}^2 + 
\|A^{1/2} f\|_{2}^2 
\, \big]^{1/2}
\end{equation*}
for every $f \in D(A^{1/2})$. Further, it 
follows for $f \in D(A)$
that 
\begin{equation*}
\|A^{1/2} f\|_{2}^2 = \braket{\,f\,|\,A f \,}_{2} = 
\braket{\,{\nabla_{w}} f \, | \, (1/\bar{p}) \,  
{\nabla_{w}} f \, }_{2,n} = s(f,f) \, \, , \, \,
\end{equation*}
where the real numbers $C_1$, $C_2$ and the sesquilinear form 
$s$ are as in the proof of 
Theorem~\ref{secondorderoperator},
and hence that 
\begin{align*}
& \min\{1,1/C_2\} \, \vvvert f \vvvert_{1}^2 
 \leq 
\|f\|_{A^{1/2}}^2  
\leq \max\{1,1/C_1\} \, \vvvert f \vvvert_{1}^2 \, \, .
\end{align*}
As a consequence, the restrictions of $\|\,\,\|_{A^{1/2}}$ and
$\vvvert \, \, \vvvert_{_{1}}$ to $D(A)$
are equivalent. Since $D(A)$ is dense in  
$(D(A^{1/2}), \|\,\|_{A^{1/2}})$, 
it follows 
for $f \in 
D(A^{1/2})$ the existence of a sequence 
$f_1,f_2,\dots$ in $D(A)$ 
such that 
\begin{equation*}
\lim_{\nu \rightarrow \infty} \|f_{\nu} - f \|_{A^{1/2}} 
= 0 \, \, .
\end{equation*}
Since the inclusion of $(D(A^{1/2}), 
\|\,\|_{A^{1/2}})$ into 
$L^{2}_{\mathbb{C}}({\Omega})$ is continuous, 
this implies also that 
\begin{equation*}
\lim_{\nu \rightarrow \infty} \|f_{\nu} - f \|_{2} = 0 \, \, .
\end{equation*}
Since the restrictions of $\|\,\,\|_{A^{1/2}}$ and
$\vvvert \, \, \vvvert_{_{1}}$ to $D(A)$ 
are equivalent, it follows that $f_1,f_2,\dots$ is a Cauchy sequence
in $W^{1}_{0,{\mathbb{C}}}({\Omega})$ and hence convergent to some 
$\bar{f} \in W^{1}_{0}({\Omega})$. Since the embedding of 
$(W^{1}_{\mathbb{C}}({\Omega}),\vvvert \, \, \vvvert_{1})$ into 
$L^{2}_{\mathbb{C}}({\Omega})$ is continuous, it follows also that 
\begin{equation*}
\lim_{\nu \rightarrow \infty} \|f_{\nu} - \bar{f} \, \|_{2} = 0 
\end{equation*}
and hence that $f = \bar{f} \in  W^{1}_{0}({\Omega})$. Further, 
it follows that 
\begin{align*}
& \min\{1,1/C_2\} \, \vvvert f \vvvert_{1}^2 
 \leq 
\|f\|_{A^{1/2}}^2  
\leq \max\{1,1/C_1\} \, \vvvert f \vvvert_{1}^2 
\end{align*}
and hence that $\|\,\,\|_{A^{1/2}}$ and the restriction 
of $\vvvert \, \, \vvvert_{_{1}}$ to $D(A^{1/2})$ are 
equivalent. 
Since according to the 
proof of Theorem~\ref{secondorderoperator}, $D(A)$ is a dense subspace 
of $(W^1_{0,{\mathbb{C}}}({\Omega}),$ $\vvvert \,\, \vvvert_{1})$,
we conclude that $D(A^{1/2}) = W^{1}_{0,{\mathbb{C}}}({\Omega})$
and that $A^{1/2} : W^1_{0,{\mathbb{C}}}({\Omega}) \rightarrow 
L^{2}_{\mathbb{C}}({\Omega})$ is continuous. From this, we conclude that 
statement of (ii) as follows. For this, let $\Lambda \in
{\mathbb{R}} \setminus \sigma(A)$ such that $\Lambda > 
\max\{0,\lambda\}$.
Since ${\mathbb{R}} \setminus \sigma(A)$ is open, such 
$\Lambda$ exists. We note that 
$D(A)$ is a core also for $A^{1/2}_{\lambda}$ and 
$A^{1/2}_{\Lambda}$. For 
instance, this follows by Theorem~3.1.9 in \cite{beyer}. Hence 
$D(A)$ is dense in the Banach spaces 
$(D(A^{1/2}_{\lambda}), \|\,\|_{A^{1/2}_{\lambda}})$, 
$(D(A^{1/2}_{\Lambda}), \|\,\|_{A^{1/2}_{\Lambda}})$,
where 
\begin{equation*}
\|f\|_{A^{1/2}_{\lambda}} := \big[ \, \|f\|_{2}^2 + 
\|A^{1/2}_{\lambda} f\|_{2}^2 
\, \big]^{1/2} \, \, , \, \, 
\|g\|_{A^{1/2}_{\Lambda}} := \big[ \, \|g\|_{2}^2 + 
\|A^{1/2}_{\Lambda} g\|_{2}^2 
\, \big]^{1/2} \, \, ,
\end{equation*}
for all $f \in D(A^{1/2}_{\lambda})$ and 
$g \in D(A^{1/2}_{\Lambda})$. Further, it follows 
for every 
$f \in D(A)$ that 
\begin{align*}
& \|f\|_{A^{1/2}_{\lambda}}^2 = 
\|A_{\lambda}^{1/2} \! f \|_{2}^{2} + 
\|f \|_{2}^{2}
= \braket{f| A_{\lambda}f}_{2}  + 
\|f \|_{2}^{2} \\
& =
\braket{f| A_{\Lambda}f}_{2} +  
\|f \|_{2}^{2} + (\Lambda - \lambda)
\|f \|_{2}^{2} = \|f\|_{A^{1/2}_{\Lambda}}^2 + 
(\Lambda - \lambda)
\|f \|_{2}^{2}
\end{align*}
and hence that 
\begin{align*}
 \|f\|_{A^{1/2}_{\lambda}}^2 \geq  \|f\|_{A^{1/2}_{\Lambda}}^2
\end{align*}
as well as that 
\begin{align*}
& \|f\|_{A^{1/2}_{\lambda}}^2 \leq  [ 1 + (\Lambda - \lambda) ] \,  
\|f\|_{A^{1/2}_{\Lambda}}^2 \, \, .
\end{align*}
Since  
$D(A)$ is dense in the Banach spaces 
$(D(A^{1/2}_{\lambda}), \|\,\|_{A^{1/2}_{\lambda}})$ and
$(D(A^{1/2}_{\Lambda}), \|\,\|_{A^{1/2}_{\Lambda}})$,
it follows that 
\begin{equation*}
D(A^{1/2}_{\lambda}) = D(A^{1/2}_{\Lambda})
\end{equation*}
as well as the equivalence of the norms $\|\,\|_{A^{1/2}_{\lambda}}$ 
and $\|\,\|_{A^{1/2}_{\Lambda}}$. In particular, this implies that 
\begin{equation*}
D(A^{1/2}_{\lambda}) = D(A^{1/2})
\end{equation*}
and the equivalence of the norms $\|\,\|_{A^{1/2}_{\lambda}}$ 
and $\|\,\|_{A^{1/2}}$. By this, the statement (ii) follows from 
the corresponding statement of (ii) for the special case 
that $\lambda = 0$. 
`(iv)': In a first step, we conclude that 
\begin{equation*}
{\nabla_{\!w}} A^{-1/2}_{\lambda} \in L(L^{2}_{\mathbb{C}}({\Omega}),
(L^{2}_{\mathbb{C}}({\Omega}))^n) \, \, .
\end{equation*}
As a consequence of the analogous properties of $A_{\lambda}$, 
$A^{1/2}_{\lambda}$ 
is 
densely-defined, linear, self-adjoint and bijective. 
Since $A^{1/2}_{\lambda}$ is self-adjoint, $A^{1/2}_{\lambda}$ is 
in particular 
closed. Further, 
according to Lemma~\ref{gradientoperatorsandadjoints}, 
the restriction of ${\nabla_{\!w}}$ to $W^{1}_{0,{\mathbb{C}}}({\Omega})$
is a closed linear operator in $L^{2}_{\mathbb{C}}({\Omega})$ with values 
in $(L^{2}_{\mathbb{C}}({\Omega}))^n$. Since $D(A^{1/2}_{\lambda}) =
W^{1}_{0,{\mathbb{C}}}({\Omega})$,
it follows from the closed graph theorem, e.g., see Theorem~3.16 
in \cite{beyer}, the existence of $C \in [0,\infty)$ such that 
\begin{equation*}
\|{\nabla_{\!w}} f \|_{2,n} \leq C \, \|A^{1/2}_{\lambda} f \|_{2}
\end{equation*}
for all $f \in D(A^{1/2}_{\lambda})$. As a consequence, it 
follows for every 
$f \in L^{2}_{\mathbb{C}}({\Omega})$ that 
\begin{equation*}
\|{\nabla_{\!w}} A^{-1/2}_{\lambda} f \|_{2,n} \leq C \, \|A_{\lambda}
^{1/2} 
A^{-1/2}_{\lambda} f \|_{2}
= C \, \|f \|_{2} 
\end{equation*}
and hence that ${\nabla_{\!w}} A^{-1/2}_{\lambda} \in L(L^{2}_{\mathbb{C}}({\Omega}),
(L^{2}_{\mathbb{C}}({\Omega}))^n)$. In a second step, we conclude that 
\begin{equation*}
\overline{A^{-1/2}_{\lambda} {\nabla_{0}}^{*}} \in 
L((L^{2}_{\mathbb{C}}({\Omega}))^n,L^{2}_{\mathbb{C}}({\Omega})) \, \, .
\end{equation*}
$A^{-1/2}_{\lambda} {\nabla_{0}}^{*}$ is a densely-defined, linear
operator in $(L^{2}_{\mathbb{C}}({\Omega}))^n$ with values in $L^{2}_{\mathbb{C}}({\Omega})$.
Further, it follows for 
$f \in L^{2}_{\mathbb{C}}({\Omega})$ that 
\begin{equation*}
\braket{\, f \, | \, A^{-1/2}_{\lambda} {\nabla_{0}}^{*} q\,}_{2} = 
\braket{\,  A^{-1/2}_{\lambda} \! f \, | \, {\nabla_{0}}^{*} q\,}_{2} =
\braket{\,{\nabla_{\!w}} A^{-1/2}_{\lambda} f\, | \, q \, }_{2,n} 
\end{equation*}
for every $q \in D({\nabla_{0}}^{*})$ and hence that 
$f \in D((A^{-1/2}_{\lambda} {\nabla_{0}}^{*})^{*})$ as well as that 
\begin{equation*}
(A^{-1/2}_{\lambda} {\nabla_{0}}^{*})^{*} f = {\nabla_{\!w}} 
A^{-1/2}_{\lambda} f \, \, .
\end{equation*} 
As a consequence, 
\begin{equation*}
(A^{-1/2}_{\lambda} {\nabla_{0}}^{*})^{*} = {\nabla_{\!w}} 
A^{-1/2}_{\lambda} \, \, .
\end{equation*}
In particular, $A^{-1/2}_{\lambda} {\nabla_{0}}^{*}$ is closable and 
\begin{equation*}
\overline{A^{-1/2}_{\lambda} {\nabla_{0}}^{*}} = 
(\,{\nabla_{\!w}} A^{-1/2}_{\lambda}\,)^{*} \in  L((L^{2}_{\mathbb{C}}({\Omega}))^{n},
L^{2}_{\mathbb{C}}({\Omega})) \, \, . 
\end{equation*}
Further, we note that
\begin{equation*}
{\nabla_{\!w}} A^{-1}_{\lambda} {\nabla_{0}}^{*} f = 
{\nabla_{\!w}} A^{-1/2}_{\lambda}  A^{-1/2}_{\lambda} {\nabla_{0}}^{*} f =
{\nabla_{\!w}} A^{-1/2}_{\lambda} \, (\,{\nabla_{\!w}} 
A^{-1/2}_{\lambda}\,)^{*} f
\end{equation*} 
for every $f \in D({\nabla_{0}}^{*})$. Hence it 
follows that $\overline{{\nabla_{\!w}} A^{-1}_{\lambda} {\nabla_{0}}^{*}
}$ is a positive self-adjoint element of $L((L^{2}_{\mathbb{C}}({\Omega}))^{n},
(L^{2}_{\mathbb{C}}({\Omega}))^n)$\,.
\end{proof}

By help of the previous lemma, we can now estimate
the size of the gap around $0$ in the spectrum of $A$ and 
give a representation of the resolvent 
of $\hat{A}$ in terms of the resolvent of $A$, i.e., 
(\ref{resolventCharacterization}).
 
\begin{theorem} \label{Ahatbijective}
Let $\bar{p} \in {\cal L}$, $C_1, C_2 \in {\mathbb{R}}$ 
satisfy $C_2 \geq C_1 > 0$ and be such that   
$C_1 \leq \bar{p} \leq C_2$ a.e. on $\Omega$. Further,
let $j \in \{1,\dots,n \}$ be such that ${\bf n}_{j} \neq 0$.
Then the interval 
\begin{equation*}
J := (\,-C_1 \, , \, c^2/(c + C_2)
\, )
\end{equation*}
is contained in the resolvent set of $\hat{A}$. In particular
for $\lambda \in J$,  
$(\hat{A} - \lambda)^{-1}$ is given by 
\begin{equation} \label{resolventCharacterization}
\begin{array}{l}
(\hat{A} - \lambda)^{-1}(f,g) =  \\
\left( (A_{{\bar p}_{\lambda}} - \lambda)^{-1} f + 
\overline{(A_{{\bar p}_{\lambda}} - \lambda)^{-1} 
{\nabla_{\!0}}^{*}} \, p_{\lambda} g \, , \, \right. \\
\left. 
- p_{\lambda} g + p_{\lambda} {\nabla_{w}} (A_{{\bar p}_{\lambda}} - \lambda)^{-1} 
f +  p_{\lambda} 
\overline{{\nabla_{w}} (A_{{\bar p}_{\lambda}} - \lambda)^{-1} 
{\nabla_{\!0}}^{*}} \, p_{\lambda} g \right)
\end{array}
\end{equation}
for all $(f,g) \in L^2_{\mathbb{C}}(\Omega) \times 
(L^2_{\mathbb{C}}(\Omega))^{n}$, where 
$\bar{p}_{\lambda} := \bar{p} + \lambda, p_{\lambda} := 
1/{\bar p}_{\lambda}$, and  
$A_{{\bar p}_{\lambda}}$ is the operator corresponding to 
$\bar{p}_{\lambda}$ according to Definition~\ref{defA}.
\end{theorem}

\begin{proof}
For this, let $\lambda \in J$.
Then,
\begin{equation*} 
0 < \lambda + C_1 
\leq \bar{p} + \lambda \leq \lambda + C_2 
\end{equation*}
and $\bar{p}_{\lambda} := \bar{p} + \lambda \in {\cal L}$. 
Further, we denote 
by $A_{{\bar p}_{\lambda}}$ the operator corresponding to 
$\bar{p}_{\lambda}$ according to Definition~\ref{defA}.
As a consequence of Theorem~\ref{Abijective}, the spectrum of 
$A_{{\bar p}_{\lambda}} - \lambda$ is contained in the interval 
\begin{equation*} 
\big[ c^2 (\lambda + 
C_{2})^{-1} - \lambda,\infty \big) \, \, ,
\end{equation*}
The inequality
\begin{equation*}
c^2 (\lambda + C_{2})^{-1} - \lambda > 0 
\end{equation*}
is equivalent to 
\begin{equation*}
\left(\lambda + \frac{C_{2}}{2} \right)^2 - \frac{C_{2}^2}{4} =
\lambda (\lambda + C_{2}) < c^2  \, \, .
\end{equation*}
The last is equivalent to 
\begin{align*}
-\sqrt{c^2 + \frac{C_{2}^2}{4}} - \frac{C_2}{2} 
< \lambda < \sqrt{c^2 + \frac{C_{2}^2}{4}} 
- \frac{C_2}{2} \, \, .
\end{align*}
We note that 
\begin{align*}
& \sqrt{c^2 + \frac{C_{2}^2}{4}} - \frac{C_2}{2} = 
\frac{c^2}{\sqrt{c^2 + \frac{C_{2}^2}{4}} + \frac{C_2}{2}} 
\geq \frac{c^2}{c+ C_2}
\end{align*}
and that 
\begin{equation*}
-\sqrt{c^2 + \frac{C_{2}^2}{4}} - \frac{C_2}{2} \leq  - C_2 \leq 
- C_1 \, \, . 
\end{equation*}
Hence it  follows that $A_{{\bar p}_{\lambda}} - \lambda$
is self-adjoint, strictly positive and  
bijective. We define 
the bounded linear operator 
$B \in L(L^{2}_{\mathbb{C}}({\Omega}) \times 
(L^{2}_{\mathbb{C}}({\Omega}))^{n},L^{2}_{\mathbb{C}}({\Omega}) \times (L^{2}_{\mathbb{C}}({\Omega}))^n)$
by 
\begin{align*}
& B (f,g) := 
((A_{{\bar p}_{\lambda}} - \lambda)^{-1} f + 
\overline{(A_{{\bar p}_{\lambda}} - \lambda)^{-1} 
{\nabla_{\!0}}^{*}} \, p_{\lambda} g \, , \, \\
& 
- p_{\lambda} g + p_{\lambda} {\nabla_{w}} (A_{{\bar p}_{\lambda}} - \lambda)^{-1} 
f +  p_{\lambda} 
\overline{{\nabla_{w}} (A_{{\bar p}_{\lambda}} - \lambda)^{-1} 
{\nabla_{\!0}}^{*}} \, p_{\lambda} g)
\end{align*}
for all $(f,g) \in L^2_{\mathbb{C}}(\Omega) \times 
(L^2_{\mathbb{C}}(\Omega))^{n}$, where $ p_{\lambda} := 
1/{\bar p}_{\lambda} \in L^{\infty}(\Omega)$. 
Further, we define the subspace 
$D$ of $L^2_{\mathbb{C}}(\Omega) \times (L^2_{\mathbb{C}}(\Omega))^{n}$
by 
\begin{equation*}
D := \{(f,g) \in L^{2}_{\mathbb{C}}({\Omega}) \times 
(L^{2}_{\mathbb{C}}({\Omega}))^{n} : p_{\lambda} 
g \in D({\nabla_{\!0}}^{*}) \} \, \, .
\end{equation*}
We note that the subspace  
\begin{equation*}
\{\bar{p}_{\lambda} g : g \in C_{0}^{\infty}(\Omega,{\mathbb{C}})\}
\end{equation*}
of $L^2_{\mathbb{C}}(\Omega)$ is dense in $L^2_{\mathbb{C}}(\Omega)$.
For the proof, let $f \in L^2_{\mathbb{C}}(\Omega)$. Since 
$p_{\lambda} \in L^{\infty}(\Omega)$, 
$p_{\lambda} f \in L^2_{\mathbb{C}}(\Omega)$. Further, since 
$C_{0}^{\infty}(\Omega,{\mathbb{C}})$ is dense 
in $L^2_{\mathbb{C}}(\Omega)$, there exists a sequence 
$f_1,f_2,\dots$ in $C_{0}^{\infty}(\Omega,{\mathbb{C}})$ such that 
\begin{equation*}
\lim_{\nu \rightarrow \infty} \|f_{\nu} - p_{\lambda} 
f\|_{2} = 0 \, \, .
\end{equation*} 
Since for every $\nu \in {\mathbb{N}}^{*}$
\begin{equation*}
\|\bar{p}_{\lambda} f_{\nu} - f\|_{2} = 
\|\bar{p}_{\lambda} (f_{\nu} - p_{\lambda} f) \|_{2} \leq 
\|\bar{p}_{\lambda}\|_{\infty} \|f_{\nu} - p_{\lambda} 
f\|_{2} \, \, , 
\end{equation*}
it follows that 
\begin{equation*}
\lim_{\nu \rightarrow \infty} \|\bar{p}_{\lambda} 
f_{\nu} - f\|_{2} = 0 \, \, .
\end{equation*}
Hence it follows also that $D$ is dense in 
$L^2_{\mathbb{C}}(\Omega) \times (L^2_{\mathbb{C}}(\Omega))^{n}$.
Further, for $(f,g) \in D$, it follows that
$B(f,g) \in D(\hat{A})$ and that 
\begin{align*}
& {\nabla_{\!0}}^{*} 
[\, - p_{\lambda} g + p_{\lambda} {\nabla_{w}} (A_{{\bar p}_{\lambda}} - \lambda)^{-1} 
f +  p_{\lambda} 
\overline{{\nabla_{w}} (A_{{\bar p}_{\lambda}} - \lambda)^{-1} 
{\nabla_{\!0}}^{*}} \, p_{\lambda} g \, ] \\
& - \lambda 
[ \, (A_{{\bar p}_{\lambda}} - \lambda)^{-1} f + 
\overline{(A_{{\bar p}_{\lambda}} - \lambda)^{-1} 
{\nabla_{\!0}}^{*}} \, p_{\lambda} g \, ] \\
& = - {\nabla_{\!0}}^{*} p_{\lambda} g + 
A_{{\bar p}_{\lambda}} (A_{{\bar p}_{\lambda}} - \lambda)^{-1} f
+ A_{{\bar p}_{\lambda}} (A_{{\bar p}_{\lambda}} - \lambda)^{-1}
{\nabla_{\!0}}^{*} \, p_{\lambda} g \\
& \phantom{=} \,\, - \lambda 
(A_{{\bar p}_{\lambda}} - \lambda)^{-1} f - \lambda
(A_{{\bar p}_{\lambda}} - \lambda)^{-1} 
{\nabla_{\!0}}^{*} \, p_{\lambda} g \, ] = f 
\end{align*}
and that 
\begin{align*}
& {\nabla_{\!w}} \, [ \, (A_{{\bar p}_{\lambda}} - \lambda)^{-1} f + 
\overline{(A_{{\bar p}_{\lambda}} - \lambda)^{-1} 
{\nabla_{\!0}}^{*}} \, p_{\lambda} g \, ] \\
& - 
{\bar{p}}_{\lambda} [\, - p_{\lambda} g + p_{\lambda} {\nabla_{w}} (A_{{\bar p}_{\lambda}} - \lambda)^{-1} 
f +  p_{\lambda} 
\overline{{\nabla_{w}} (A_{{\bar p}_{\lambda}} - \lambda)^{-1} 
{\nabla_{\!0}}^{*}} \, p_{\lambda} g \, ] \\
& = 
 {\nabla_{\!w}} (A_{{\bar p}_{\lambda}} - \lambda)^{-1} f + 
{\nabla_{\!w}} (A_{{\bar p}_{\lambda}} - \lambda)^{-1} 
{\nabla_{\!0}}^{*} \, p_{\lambda} g  \\
& \phantom{=} \,\, \,\,   
g - {\nabla_{w}} (A_{{\bar p}_{\lambda}} - \lambda)^{-1} 
f - {\nabla_{w}} (A_{{\bar p}_{\lambda}} - \lambda)^{-1} 
{\nabla_{\!0}}^{*} \, p_{\lambda} g  = g \, \, . 
\end{align*}
Hence it follows that   
\begin{equation} 
(\hat{A} - \lambda) B (f,g) = (f,g) \, \, .
\end{equation}
Further, since $D$ is dense in 
$L^2_{\mathbb{C}}(\Omega) \times (L^2_{\mathbb{C}}(\Omega))^{n}$,
for $(f,g) \in 
L^2_{\mathbb{C}}(\Omega) \times (L^2_{\mathbb{C}}(\Omega))^{n}$,
there is a sequence $(f_1,g_1),(f_2,g_2),\dots$ in $D$ that is 
convergent to $(f,g)$. Since $B$ is bounded, the corresponding 
sequence $B(f_1,g_1),B(f_2,g_2),\dots$ is convergent to $B(f,g)$. 
Since $\hat{A}$ is in particular closed, it follows that 
$B(f,g) \in D(\hat{A})$
and that 
\begin{equation}
(\hat{A} - \lambda) B (f,g) = (f,g)  \, \, . 
\end{equation} 
Therefore, $\hat{A} - \lambda$ is surjective and hence also 
bijective.  
\end{proof}

By help of the previous theorem, the next result follows 
by application of a well-known criterion 
for the strong resolvent convergence of sequences of self-adjoint 
operators.   

\begin{theorem} \label{stongconvergence}
Let ${\bar{p}}_{1}, {\bar{p}}_{2}, \dots $ be a uniformly 
bounded sequence 
in ${\cal L}$ for which there is $\varepsilon >0$ such 
that ${\bar{p}}_{\nu} \geq \varepsilon$ for 
all $\nu \in {\mathbb{N}}^{*}$
and 
which converges a.e. pointwise on 
$\Omega$ to ${\bar p}_{\infty} \in {\cal L}$. In addition, let 
${\hat{A}}_1,{\hat{A}}_2,\dots$ be the associated 
sequence of self-adjoint
linear operators and ${\hat{A}}_{\infty}$ be 
the self-adjoint linear operator associated to 
${\bar{p}_{\infty}}$. Then  
\begin{equation*}
s\!-\!\lim_{\nu \rightarrow \infty} {\hat{A}}_{\nu}^{-1} 
= {\hat{A}}_{\infty}^{-1} \, \, . 
\end{equation*}
\end{theorem}

\begin{proof}
By application of Lebesgue's dominated convergence theorem, 
it follows that 
\begin{equation*}
\lim_{\nu \rightarrow \infty} \|{\hat{A}}_{\nu} (u,q) - 
{\hat{A}}_{\infty}(u,q)\| = 0 \, \, , 
\end{equation*}
for all $(u,q) \in L^2_{\mathbb{C}}(\Omega) \times 
(L^2_{\mathbb{C}}(\Omega))^{n}$, 
where $\|\,\|$ denotes the norm on $L^2_{\mathbb{C}}(\Omega) \times 
(L^2_{\mathbb{C}}(\Omega))^{n}$. From this, the statement follows 
from a well-known criterion for strong resolvent convergence of a 
sequence of self-adjoint linear operators, e.g., see part (i)
of Theorem~9.16 in \cite{weidmann}.
\end{proof}

\begin{corollary} \label{stongconvergencecor}
Let ${\bar{p}}_{1}, {\bar{p}}_{2}, \dots $ be a uniformly 
bounded sequence 
in ${\cal L}$  for which there is $\varepsilon >0$ such 
that ${\bar{p}}_{\nu} \geq \varepsilon$ for 
all $\nu \in {\mathbb{N}}^{*}$
and which converges a.e. pointwise on 
$\Omega$ to $\bar{p}_{\infty} \in {\cal L}$. In addition, let 
$A_1,A_2,\dots$ be the associated 
sequence of self-adjoint
linear operators and $A_{\infty}$ be 
the self-adjoint linear operator associated to 
${\bar{p}_{\infty}}$. Finally, let $f \in L^2_{\mathbb{C}}(\Omega)$. Then
\begin{equation} \label{mainresult}
\lim_{\nu \rightarrow \infty}
\|A_{\nu}^{-1} f - A_{\infty}^{-1} f \|_{2} = 
\lim_{\nu \rightarrow \infty} 
\|p_{\nu} \nabla_{\! w} A_{\nu}^{-1} f - p_{\infty} 
\nabla_{\! w} A_{\infty}^{-1} f \|_{2,n} = 0 \, \, ,
\end{equation}
where $p_{\nu} := 1 / {\bar{p}}_{\nu}$ for every $\nu \in 
{\mathbb{N}}^{*}$ and $p_{\infty} := 1 / {\bar{p}}_{\infty}$.
\end{corollary}

\begin{proof}
By Theorem~\ref{stongconvergence}, it follows that 
\begin{equation*}
\lim_{\nu \rightarrow \infty} {\hat{A}}_{\nu}^{-1}(f,0) 
= {\hat{A}}_{\infty}^{-1}(f,0)  \, \, , 
\end{equation*}
where 
${\hat{A}}_1,{\hat{A}}_2,\dots$ is the associated 
sequence of self-adjoint
linear operators to ${\bar{p}}_{1}, {\bar{p}}_{2}, \dots $
and ${\hat{A}}_{\infty}$ is the self-adjoint linear operator 
associated to ${\bar{p}}$. Hence (\ref{mainresult}) follows 
by Theorem~\ref{Ahatbijective}.
\end{proof}

\section{The one-dimensional case}

In the special cases that $\Omega$ is given by a non-empty bounded
open interval of ${\mathbb{R}}$, ${\hat{A}}^{-1}$ can be explicitly
calculated. This is somewhat surprising since in this case the
corresponding $A$ is a Sturm-Liouville operator and the standard
method of calculating its inverse, e.g., see Theorem~8.26
in \cite{weidmann}, seems not directly applicable for general $p \in
{\cal L}$. Still, by a direct calculation of ${\hat{A}}^{-1}$,
one can give an explicit expression of $A^{-1}$ by using 
(\ref{resolventCharacterization}).

\begin{theorem} \label{Ahatbijectivein1dim}
Let $a,b \in {\mathbb{R}}$ such that $a < b$ and $\Omega := I := 
(a,b)$. Further, let $\bar{p} \in L^{\infty}({\Omega}) \setminus \{0\}$ 
be a.e. positive. Then 
$\hat{A}$ is bijective and has a purely discrete spectrum. 
In particular, ${\hat{A}}^{-1}$ is given by 
\begin{equation*}
(u,q) := {\hat{A}}^{-1}(f,g)
\end{equation*}
for every $(f,g) \in (L^2_{\mathbb{C}}(I))^2$, where 
\begin{align*}
& u(x) = \int_{a}^{x} g(y) \, dy + \int_{a}^{x} 
\left[ \int_{y}^{x} \bar{p}(u) \, du  \right] \! f(y) \, dy \\
& \phantom{u(x) =} + 
\|\bar{p}\|_{1}^{-1} \int_{a}^{x} \bar{p}(y) \, dy \, \left\{
\int_{a}^{b} \left[ \int_{y}^{b} \bar{p}(x) \, dx  \right] \! f(y) \,
dy - \int_{a}^{b} g(y) \, dy \right\} \\
& q(x) = - \int_{a}^{x} f(y) \, dy + \|\bar{p}\|_{1}^{-1}  \left\{
\int_{a}^{b} \left[ \int_{y}^{b} \bar{p}(x) \, dx  \right] \! f(y) \,
dy - \int_{a}^{b} g(y) \, dy \right\} 
\end{align*}
for every $x \in I$. Also, ${\hat{A}}^{-1}$ satisfies
\begin{equation*}
\|{\hat{A}}^{-1}\| \leq 2 \, (b - a) \, \|\bar{p}\|_{1}^{-1} 
(\, 1 + \|\bar{p}\|_{1})^2  \, \, .
\end{equation*}
\end{theorem}

\begin{proof}
For this, we define the derivative operator
\begin{equation*}
D_{I} : C^{\, \infty}_{0}(I, \, {\mathbb{C}}) 
\rightarrow L^2_{\mathbb{C}}(I)
\end{equation*} 
by 
$D_{I} f := f^{\, \prime}$ 
for every $f \in  C^{\, \infty}_{0}(I, \, {\mathbb{C}})$. 
In a first step, we prove an auxiliary result. For this, 
let $f \in L^{2}_{\mathbb{C}}(I)$ and $h \in 
C(\bar{I},{\mathbb{C}})$ be defined by 
\begin{equation*}
h(x) := \int_{a}^{x} f(y) \, dy 
\end{equation*} 
for every $x \in I$.
Further, let $\varphi \in C_{0}^{\infty}(I,{\mathbb{C}})$. By 
Fubini's theorem and change of variables, it follows that 
\begin{align*}
& \braket{\, h \, | \, D_{I} \varphi \, }_{2} =
\braket{\, h \, | \, \varphi^{\, \prime} \, }_{2} =
\int_{a}^{b} h^{*}(x) \, \varphi^{\, \prime}(x) \, dx  =
\int_{a}^{b} \left[ \, \int_{a}^{x} \varphi^{\, \prime}(x) 
f^{*}(y) \, dy \right] \! dx \\
& = \int_{\{(x,y) \in {\mathbb{R}}^2 : a \leq x \leq b \wedge 
a \leq y \leq x\}} \varphi^{\, \prime}(x) 
f^{*}(y) \, dx dy \\
& = \int_{\{(x,y) \in {\mathbb{R}}^2 : a \leq y \leq b \wedge 
y \leq x \leq b \}} \varphi^{\, \prime}(x) 
f^{*}(y) \, dx dy \\
& = \int_{a}^{b} \left[ \, \int_{y}^{b} \varphi^{\, \prime}(x) 
f^{*}(y) \, dx \right] dy = 
- \int_{a}^{b} \varphi(y) f^{*}(y) \, dy = -
\braket{\, f \, | \, \varphi \,}_{2}
\end{align*}
and hence that  
\begin{equation*} 
h \in W^{1}_{\mathbb{C}}(I) \, \, \textrm{and} \, \, 
D_{I}^{*} h =  - f \, \, .
\end{equation*}
With the help of the previous auxiliary result, we proceed
in the proof of the lemma. 
For this, we define for every $(f,g) \in (L^{2}_{\mathbb{C}}(I))^2$,
a corresponding $B(f,g) = (u,q)$ by 
\begin{equation*}
q(x) := q_{0}(x) + c \, \, , \, \, 
u(x) := \int_{a}^{x} \left[ \,   
g(y) + \bar{p}(y) (q_{0}(y) + c) 
\right] dy \, \, , 
\end{equation*}
where 
\begin{equation*}
q_{0}(x) := - \int_{a}^{x} f(y) \, dy \, \, , \, \,
c := - \|\bar{p}\|_{1}^{-1} \,  
\int_{a}^{b}  \left[ \,  
g(y) + \bar{p}(y) q_{0}(y) \, \right] dy 
\end{equation*}
for every $x \in I$. By help of the auxiliary result above, it follows 
that $(u,q) \in (W^1_{\mathbb{C}}(I) \cap 
C(\bar{I},{\mathbb{C}})) \times D(D^{*}_{I})$ and that 
\begin{align*}
D^{*}_{I} q = f \, \, , \, \, 
- D^{*}_{I} u - \bar{p} q = g + \bar{p} q - \bar{p} q = g \, \, .
\end{align*}
In addition,
\begin{align*}
& u_{b} = \int_{a}^{b} \left[ \,  
g(y) + \bar{p}(y) (q_{0}(y) + c) \, 
\right] dy \\
& = \int_{a}^{b} \left( \,  
g(y) + \bar{p}(y) q_{0}(y) \, 
\right) dy + c \int_{a}^{b} \bar{p}(y) \, dy = 0 \, \, . 
\end{align*}
As a consequence, 
\begin{equation*}
u_{a} = u_{b} = 0 \, \, .
\end{equation*}
From the last, it follows also that $u \in W^1_{0, \mathbb{C}}(I)$. 
For the proof, let $\varphi \in C^{\infty}({\mathbb{R}})$ be 
such that 
\begin{equation*}
\varphi (\,(-\infty,0]\,) \subset \{0\} \, \, , \, \,
\varphi (\,[1,\infty)\,) \subset \{1\} \, \,  , \, \, 
\textrm{Ran} \, \varphi  \subset [0,1] \, \, .
\end{equation*}
Such a function is easy to construct. 
For $\nu \in {\mathbb{N}}^{*}$, 
we define $\varphi_{\nu} 
\in C^{\infty}({\mathbb{R}})$ by 
\begin{equation*}
\varphi_{\nu}(x) := \varphi(\, \nu (x-a)
(b - x) - 1 \,)
\end{equation*}
for every $x \in I$. 
Then it follows for every 
$\nu \in {\mathbb{N}}^{*}$ satisfying 
$\nu \geq b-a$
and $x \in I$ that 
\begin{equation*}
\varphi_{\nu}(x) = \begin{cases}
0 & 
\text{if $ x \in I \, \setminus \, 
(a+\nu^{-2},b-\nu^{-2})$}  \\
1
& \text{if $(x-a)(b-x) \geq 2 \,\nu^{-1}$}
\end{cases}
\end{equation*}
and hence that $\varphi_{\nu} \in  C^{\infty}_{0}(I,{\mathbb{R}})$
as well as $\textrm{Ran}(\varphi_{\nu}) \subset [0,1]$. In particular, 
\begin{align*}
& |\,(\,x - a\,)( \, b - x \, ) 
\, {\varphi}_{\nu}^{\, \prime}(x)\,| \\
& \leq 3 \nu \, (a+b) \, 
(\,x - a\,)( \, b - x \, ) \, 
 \cdot 
|\,{\varphi}^{\, \prime}(\,\nu (x-a)
(b - x) - 1 \,) \, | \\
& \leq 3 \nu \, (a+b) \, 
(\,x - a\,)( \, b - x \, ) \cdot \|{\varphi}^{\, \prime} 
\|_{\infty} \cdot 
\chi_{_{\,\{x \in I: (x-a)(b-x) \leq 2 / \nu \}\,}}(x) \\
& \leq 6 \, (a+b) \cdot \|{\varphi}^{\, \prime} 
\|_{\infty} \cdot 
\chi_{_{\,\{x \in I: (x-a)(b-x) \leq 2 / \nu \}\,}}(x)
\end{align*}
for all $x \in I$. An application of 
Lebesgue's 
dominated convergence theorem leads to  
\begin{equation*}
\lim_{\nu \rightarrow \infty} \|{\varphi}_{\nu}u - u\|_{2} = 
\lim_{\nu \rightarrow \infty} \| 
D_{I}^{*} {\varphi}_{\nu}u - 
D_{I}^{*}u \|_{2} =
\lim_{\nu \rightarrow \infty} \| 
{\varphi}_{\nu} D_{I}^{*} u + {\varphi}_{\nu}^{\, \prime}
u - D_{I}^{*}u \|_{2} = 0 \, \, .
\end{equation*}
Hence it follows that $u \in  W^1_{0,{\mathbb{C}}}(I)$ and further 
that $(u,q) \in D(\hat{A})$ and 
\begin{equation*}
\hat{A} B (f,g) = \hat{A}(u,q) = (f,g) \, \, . 
\end{equation*}
Further, we conclude by Fubini's theorem that  
\begin{align*}
& c = - \|\bar{p}\|_{1}^{-1}  \! 
\int_{a}^{b} g(y) \, dy + \|\bar{p}\|_{1}^{-1}  \! 
\int_{a}^{b} \bar{p}(x) \left[ \int_{a}^{x} f(y) \, dy  \right]
dx \\
& = \|\bar{p}\|_{1}^{-1}  \left\{
\int_{a}^{b} \left[ \int_{y}^{b} \bar{p}(x) \, dx  \right] \! f(y) \,
dy - \int_{a}^{b} g(y) \, dy \right\} \, \, .
\end{align*}
This implies that 
\begin{equation*}
q(x) = - \int_{a}^{x} f(y) \, dy + \|\bar{p}\|_{1}^{-1}  \left\{
\int_{a}^{b} \left[ \int_{y}^{b} \bar{p}(x) \, dx  \right] \! f(y) \,
dy - \int_{a}^{b} g(y) \, dy \right\} 
\end{equation*}
for every $x \in I$.
Further, again by Fubini's theorem, it follows that 
\begin{align*}
& u(x) = \int_{a}^{x} g(y) \, dy 
+ \int_{a}^{x} \bar{p}(y) \, q_{0}(y) \, 
dy + c \int_{a}^{x} \bar{p}(y) \, dy \\
& = \int_{a}^{x} g(y) \, dy 
+ \int_{a}^{x} \bar{p}(u) \left[ \int_{a}^{u} f(y) \, dy  \right]
du + c \int_{a}^{x} \bar{p}(y) \, dy \\
& = \int_{a}^{x} g(y) \, dy + \int_{a}^{x} 
\left[ \int_{y}^{x} \bar{p}(u) \, du  \right] \! f(y) \, dy \\
& \phantom{=} \, \, + 
\|\bar{p}\|_{1}^{-1} \int_{a}^{x} \bar{p}(y) \, dy \, \left\{
\int_{a}^{b} \left[ \int_{y}^{b} \bar{p}(x) \, dx  \right] \! f(y) \,
dy - \int_{a}^{b} g(y) \, dy \right\}
\end{align*}
for every $x \in I$. In addition, by Hoelder's inequality, we conclude that 
\begin{align*}
& |u(x)| \leq 2 \, (b - a)^{1/2} \left[ \, 
\|\bar{p}\|_{1}  \, \|f\|_{2} + \|g\|_{2} \, \right] \leq 
2 \, (b - a)^{1/2} \, (\,1 + \|\bar{p}\|_{1}) \, \|(f,g)\|
\\
& |q(x)| \leq (b - a)^{1/2} \left[ \, 2 \, \|f\|_{2} +  
\|\bar{p}\|_{1}^{-1}
\|g\|_{2} \, \right] \\
& \phantom{|q(x)|} \leq 2 \, (b - a)^{1/2} \, \|\bar{p}\|_{1}^{-1} \, 
(\,1 + \|\bar{p}\|_{1}) \, 
\|(f,g)\|_{2}
\end{align*}
for every $x \in I$. The last implies  
\begin{equation*}
\|u\|_{2} \leq 2 \, (b - a) \, (\,1 + \|\bar{p}\|_{1}) \, \|(f,g)\|
\, \, , \, \, \|q\|_{2} \leq 2 \, (b - a) \, \|\bar{p}\|_{1}^{-1} \, 
(\,1 + \|\bar{p}\|_{1}) \, 
\|(f,g)\|_{2}
\end{equation*}
and 
\begin{equation*}
\|(u,q)\| \leq 2 \, (b - a) \, \|\bar{p}\|_{1}^{-1} 
(\,1 + \|\bar{p}\|_{1})^2 \, 
\|(f,g)\|_{2} \, \, .
\end{equation*}
As consequence, by $((L^{2}_{\mathbb{C}}(I))^2
\rightarrow (L^2_{\mathbb{C}}(I))^2, 
(f,g) \rightarrow B(f,g))$, there is defined a compact bounded linear 
operator $B$. Since 
\begin{equation*}
\hat{A} B(f,g) = (f,g)
\end{equation*}
for every $(f,g) \in (L^{2}_{\mathbb{C}}(I))^2$, the bijectivity of
$\hat{A}$ follows as well as that ${\hat{A}}^{-1} = B$. Finally, since
${\hat{A}}^{-1}$ is compact, $\hat{A}$ has a purely discrete spectrum.
\end{proof}

\begin{corollary} 
Let $a,b,\bar{p}$ as in Theorem~\ref{Ahatbijectivein1dim}. 
Then $U_{r}(0)$, where 
\begin{equation*}
r := 2^{-1} (b-a)^{-1} \, \|\bar{p}\|_{1} 
\, (\,1 + \|\bar{p}\|_{1})^{-2} \, \, , 
\end{equation*}
is contained in the resolvent set of $\hat{A}$.
\end{corollary}

\begin{proof}
For this, let $\lambda \in U_{r}(0)$. Then 
\begin{equation*}
\hat{A} - \lambda = (1 - \lambda \, {\hat{A}}^{-1}) \hat{A} 
\, \, . 
\end{equation*}
By help of the previous Theorem~\ref{Ahatbijectivein1dim},
it follows that  
\begin{equation*}
| \lambda | \cdot \|{\hat{A}}^{-1} \| \leq 
| \lambda | / r < 1
\end{equation*}
and hence that $\hat{A} - \lambda$ is bijective. 
\end{proof}

\begin{theorem} \label{speedofconv}
Let $a,b \in {\mathbb{R}}$ such that $a < b$ and $\Omega := I := 
(a,b)$. Further, let ${\bar{p}}_{1}, {\bar{p}}_{2} 
\in L^{\infty}({\Omega}) \setminus \{0\}$ be a.e. positive
and ${\hat{A}}_1, 
{\hat{A}}_2$ be the corresponding operators. Then
\begin{equation*}
\|{\hat{A}}_1^{-1} - {\hat{A}}_2^{-1} \| \leq 
\frac{2 (b-a)}{\| {\bar{p}}_{1}\|_{1}} 
\left(2  + \|{\bar{p}}_{1}\|_{1} + \|{\bar{p}}_{2}\|_{1} + 
\frac{1}{\|{\bar{p}}_{2}\|_{1}} \right)   
\|\, {\bar{p}}_{2} - {\bar{p}}_{1} \|_{1} \, \, .
\end{equation*}
\end{theorem}

\begin{proof}
Proceeds by direct calculation. 
\end{proof}

\begin{corollary} \label{uniformconvergence1dim}
Let $a,b \in {\mathbb{R}}$ such that $a < b$ and $\Omega := I := 
(a,b)$. Further, let ${\bar{p}}_{\infty} 
\in L^{\infty}({\Omega}) \setminus \{0\}$ 
be a.e. positive.
Let ${\bar{p}}_{1}, {\bar{p}}_{2}, \dots $ be  
a sequence of a.e. positive elements 
of $L^{\infty}({\Omega}) \setminus \{0\}$ such 
that 
\begin{equation*}
\lim_{\nu \rightarrow \infty} \|{\bar{p}}_{\nu} - 
{\bar{p}}_{\infty} \|_{1} = 0 \, \, . 
\end{equation*} 
In addition, let 
${\hat{A}}_1,{\hat{A}}_2,\dots$ be the associated 
sequence of self-adjoint
linear operators and ${\hat{A}}_{\infty}$ be 
the self-adjoint linear operator associated to 
${\bar{p}_{\infty}}$. Then 
\begin{equation*}
\lim_{\nu \rightarrow \infty} \|{\hat{A}}_{\nu}^{-1} 
- {\hat{A}}_{\infty}^{-1}\| = 0  \, \, . 
\end{equation*}
\end{corollary}

\begin{proof}
The statement is a simple consequence of 
Theorem~\ref{Ahatbijectivein1dim}.
\end{proof}

\begin{corollary}  \label{uniformconvergence1dim2}
Let $a,b \in {\mathbb{R}}$ such that $a < b$ and $\Omega := I := 
(a,b)$ and $f \in L^2_{\mathbb{C}}(I)$. Further, let 
${\bar{p}}_{\infty}  \in L^{\infty}({\Omega}) \setminus \{0\}$ 
be a.e. positive and 
${\bar{p}}_{1}, {\bar{p}}_{2}, \dots $ be a sequence in 
${\cal L}$  such that 
\begin{equation*}
\lim_{\nu \rightarrow \infty} \|{\bar{p}}_{\nu} - 
{\bar{p}}_{\infty} \|_{1} = 0 \, \, . 
\end{equation*} 
 In addition, let 
$A_1,A_2,\dots$ be the sequence of self-adjoint
linear operators that is associated to 
${\bar{p}}_{1}, {\bar{p}}_{2}, \dots \, \, \, $ and 
$p_{\nu} := 1/{\bar{p}}_{\nu}$ for $\nu \in {\mathbb{N}}^{*}$.
Then 
$A_{1}^{-1}, A_{2}^{-1},\dots$ and $- p_{1} 
D_{I}^{*} A_{1}^{-1}, - p_{2} D_{I}^{*} A_{2}^{-1}, \dots$
are convergent in $L(L^2_{\mathbb{C}}(I),L^2_{\mathbb{C}}(I))$
to $B, C \in L(L^2_{\mathbb{C}}(I),L^2_{\mathbb{C}}(I))$, respectively. 
In particular, $B$ and $C$ are given by   
\begin{align*}
& (Bf)(x) = \int_{a}^{x} 
\left[ \int_{y}^{x} {\bar{p}}_{\infty}(u) \, du  \right] \! f(y) \, dy \\
& \phantom{(Bf)(x) =} + 
\|{\bar{p}}_{\infty}\|_{1}^{-1} \int_{a}^{x} {\bar{p}}_{\infty}(y) 
\, dy \,
\int_{a}^{b} \left[ \int_{y}^{b} {\bar{p}}_{\infty}(x) \, dx  
\right] \! f(y) \,
dy \, \, , \\
& (Cf)(x) = - \int_{a}^{x} f(y) \, dy + \|{\bar{p}}_{\infty}\|_{1}^{-1} 
\int_{a}^{b} \left[ \int_{y}^{b} {\bar{p}}_{\infty}(x) \, dx  
\right] \! f(y) \,
dy   
\end{align*}
for all $x \in I$ and every $f \in L^2_{\mathbb{C}}(I)$. 
\end{corollary}

\begin{proof}
The statement is a simple consequence of Theorem~\ref{Ahatbijective}
and Theorem \ref{Ahatbijectivein1dim}.
\end{proof}

\section{Concluding remarks}
It is unclear whether results similar to those of the previous 
section can be expected to hold in dimensions greater than $1$.
According to Theorem~\ref{characterizationkerAhat}
and the subsequent example, and differently to the 
situation in one dimension, ${\hat{A}}$ is not injective 
when ${\bar{p}}$ vanishes on non-empty 
open subsets of the material. Hence there 
does not seem to be an obvious candidate for a limit of
a sequence of   
${\hat{A}}^{-1}$ that is associated to a sequence in ${\cal L}$ 
approaching such ${\bar{p}}$. Therefore, it is conceivable 
that such limits show a wider variety of phenomena than those 
in one dimension. This problem deserves further 
study.   
\newline
\linebreak
A final remark concerns the fact that it cannot be expected that 
general `elliptic regularity theorems' hold for operators $A$ 
corresponding to discontinuous diffusivities $p$ as a consequence 
of the condition that every
element $u$ from the domain of such operator 
satisfies $p \, \nabla_{w} u \in D(\nabla_{0}^{*})$. For 
instance, the source function $f$ in Example~\ref{counterex} 
is in $W^{k}_{\mathbb{C}}(I)$ for every $k \in {\mathbb{N}}$, but 
$u = A^{-1} f \notin W^{2}_{\mathbb{C}}(I)$, where $I$ is the open 
interval $(-1,1)$ of ${\mathbb{R}}$.

\section{Appendix}

In the following, proofs of the 
Lemmata~\ref{partialintegration},~\ref{gradientoperatorsandadjoints} 
from Section~\ref{prerequisites} are given.

\begin{lemma} ({\bf Partial integration})
\begin{equation*} 
\braket{f|\partial^{\,e_k} g}_2 = - \braket{\partial^{\,e_k} f|g}_2  
\end{equation*}
for all $(f,g) \in  
W^1_{0,{\mathbb{C}}}(\Omega) \times  W^1_{{\mathbb{C}}}(\Omega)$
and $k \in {\mathbb{N}}^{*}$, where $e_k$ denotes the $k$-th 
canonical unit vector of ${\mathbb{R}}^{n}$.
\end{lemma}
\begin{proof}
For this, let $k \in {\mathbb{N}}^{*}$. We define the sesquilinear 
form $s :  W^1_{0,{\mathbb{C}}}(\Omega) \times  
W^1_{{\mathbb{C}}}(\Omega) \rightarrow 
{\mathbb{C}}$ by  
\begin{equation*}
s(f,g) := \braket{f|\partial^{\,e_k} g}_2 + 
\braket{\partial^{\,e_k}f|g}_2
\end{equation*}
for all $(f,g) \in  W^1_{0,{\mathbb{C}}}(\Omega) \times  
W^1_{{\mathbb{C}}}(\Omega)$. 
By the continuity of $\partial^{\,e_k}$, it follows the continuity 
of $s$ and 
by partial integration 
that $s(f,g)= 0$ for all $f \in C_0^{\infty}(\Omega,{\mathbb{C}})$
and $f \in  C^{\infty}(\Omega,{\mathbb{C}}) 
\cap W^1_{{\mathbb{C}}}(\Omega)$. 
Since  
$C_0^{\infty}(\Omega,{\mathbb{C}})\times ((C^{\infty}(\Omega,{\mathbb{C}}) \cap 
W^1_{{\mathbb{C}}}(\Omega))$ is dense in  
$W^1_{0,{\mathbb{C}}}(\Omega) \times  W^1_{{\mathbb{C}}}(\Omega)$,
this implies the vanishing of $s$ and hence the validity of  
$(\ref{partialintegration3})$
for all $(f,g) \in  W^1_{0,{\mathbb{C}}}(\Omega) \times  
W^1_{{\mathbb{C}}}(\Omega)$. 
\end{proof}

\begin{lemma} 
({\bf Adjoints of gradient operators}) 
\begin{equation*} 
({\nabla_{\!0}}^{*})^{*} = \nabla_{\!w}
\big|_{W^{1}_{0, \mathbb{C}}(\Omega)}
\, \, , \, \,
\left( \nabla_{\!w}
\big|_{W^{1}_{0, \mathbb{C}}(\Omega)}\right)^{*} = 
{\nabla_{\!0}}^{*} \, \, .
\end{equation*}
\end{lemma}

\begin{proof}
Since ${\nabla_{\!0}}^{*}$ is densely-defined, it follows that 
\begin{equation*}
({\nabla_{\!0}}^{*})^{*} = \overline{{\nabla_{\!0}}} \, \, .
\end{equation*}
For $f \in D(\,\overline{{\nabla_{\!0}}}\,)$, there exists a  
sequence $f_1,f_2,\dots$ in $C_{0}^{\infty}({\Omega},{\mathbb{C}})$
such that 
\begin{equation*}
\lim_{\nu \rightarrow \infty} \|f_{\nu} - f\|_{2} = 0 \, \, , 
\lim_{\nu \rightarrow \infty} \|{\nabla_{\!0}}f_{\nu} - 
\overline{{\nabla_{\!0}}}f\|_{2,n} =
\lim_{\nu \rightarrow \infty} \|
{\nabla_{\!w}} f_{\nu} - \overline{{\nabla_{\!0}}}f \|_{2,n}
= 0 \, \, .
\end{equation*} 
Hence it follows that $f \in W^{1}_{0, \mathbb{C}}({\Omega})$ and 
that $\overline{{\nabla_{\!0}}}f = \nabla_{\!w} f$.
As a consequence, 
\begin{equation*}
\overline{{\nabla_{\!0}}} \subset \nabla_{\!w}
\big|_{W^{1}_{0, \mathbb{C}}(\Omega)} \, \, .
\end{equation*}
Further, for $f \in W^{1}_{0, \mathbb{C}}({\Omega})$, there is 
a sequence $f_1,f_2,\dots$ in $C_{0}^{\infty}({\Omega},{\mathbb{C}})$
such that 
\begin{equation*}
\lim_{\nu \rightarrow \infty} \|f_{\nu} - f\|_{2} = 0 \, \, ,  
\lim_{\nu \rightarrow \infty} \|\nabla_{\!0} f_{\nu} 
- \nabla_{\!w} f\|_{2}
= 0 \, \, .
\end{equation*}
Hence it follows that 
\begin{equation*}
(f,\nabla_{\!w} f) \in  G(\,\overline{{\nabla_{\!0}}}\,) \, \, .
\end{equation*} 
As a consequence, 
\begin{equation*}
\nabla_{\!w}
\big|_{W^{1}_{0, \mathbb{C}}(\Omega)} \subset \overline{{\nabla_{\!0}}} \, \, .
\end{equation*}
Finally, it follows the validity of (\ref{gradientsandadjoints}, 1). The 
validity of (\ref{gradientsandadjoints}, 2)
is a simple consequence of (\ref{gradientsandadjoints}, 1) and the 
closedness 
of ${\nabla_{\!0}}^{*}$. 
\end{proof}

\section{Acknowledgment}
B. Aksoylu would like to thank Center for Computation Technology at
Louisiana State University for generous support of the research, and
A. Knyazev for his hospitality during a visit at University of
Colorado at Denver.


\begin{thebibliography}{99}
\bibitem{adams} Adams R A, Fournier J J F 2003, 
{\em Sobolev spaces}, 2nd ed., Academic Press: Amsterdam.

\bibitem{AGKS:2007}
Aksoylu B, Graham I G, Klie H, and Scheichl R 2008, 
{\em Towards a rigorously justified algebraic preconditioner for high-contrast diffusion problems}, 
Comput. Vis. Sci., {\bf 11}, 319-331, 
\newblock{doi:10.1007/s00791-008-0105-1}.

\bibitem{ABDP:1981}
Alcouffe R E, Brandt A, Dendy J E, Painter J W 1981,
{\em The multi--grid methods for the diffusion equation with strongly 
discontinuous coefficients}, 
SIAM J. Sci. Stat. Comput., {\bf 2}, 430-454.

\bibitem{bakhvalov1}Bakhvalov N S, Knyazev A V 1990, 
{\em A new iterative algorithm for solving problems of the 
fictitious flow method for elliptic equations}, Soviet Math. Dokl. 
{\bf 41}, 481-485.

\bibitem{BMMS_fosls:2005}
Berndt M, Manteuffel T A, McCormick  S F, Starke G 2005,
{\em Analysis of first-order system least squares {(FOSLS)} for elliptic problems
  with discontinuous coefficients: {Part I}}, 
SIAM J. Numer. Anal., {\bf 43}, 386-408.

\bibitem{beyer} Beyer H R 2007,
{\em Beyond partial differential equations: A course on linear and
quasi-linear abstract hyperbolic evolution equations},
Lecture Notes in Math., {\bf 1898}, Springer: Berlin.

\bibitem{braess} Braess D 2007,  
{\em Finite Elements: Theory, Fast Solvers, and 
Applications in Solid Mechanics}, 3rd ed., Cambridge University Press: 
Cambridge. 

\bibitem{BrannickEtAl_SciDac:2006}
Brannick J, Brezina M, Falgout R, Manteuffel T, McCormick S, Ruge J,
Sheehan B, Xu J, Zikatanov L 2006,
{\em Extending the applicability of multigrid methods}, 
Journal of Physics: Conference Series, {\bf 46}, 443-452,
\newblock{SciDAC 2006}.

\bibitem{CLMM:1994_foslsFirstPaper}
Cai Z, Lazarov R D, Manteuffel T A, McCormick S F 1994,
{\em First-order system least squares for second-order partial differential
  equations: {Part I}}, SIAM J. Numer. Anal., {\bf 31}, 1785-1802.

\bibitem{engel} Engel K-J, Nagel R 2000, {\em One-parameter semigroups
for linear evolution equations}, Springer: New York.

\bibitem{kato} Kato T 1966,
{\em Perturbation theory for linear operators}, Springer: New York.

\bibitem{KlWiDr:2002}
Klawonn A, Widlund O B, Dryja M 2002,
{\em {Dual-primal FETI methods for three-dimensional elliptic problems with 
heterogeneous coefficients}},
SIAM J. Numer. Anal., {\bf 40}, 159-179.

\bibitem{knyazev1} Knyazev A V 1992, 
{\em Iterative solution of PDE with strongly varying 
coefficients: Algebraic version}, in: Beauwens R, 
De Groen P (eds) 1992, {\em Iterative methods in linear algebra}, 
Elsevier: New York, 85-89. 

\bibitem{KnWi:2003}
Knyazev A and Widlund O 2003,
{\em {Lavrentiev regularization + Ritz approximation = uniform finite element 
error estimates for differential equations with rough coefficients}}, 
Math. Comp., {\bf 72}, 17-40.

\bibitem{lions} Lions J L 1973, {\em Perturbations 
Singulieres Dans Les Problemes Aux Limites Et En Controle Optimal}, 
Lecture Notes in Math., {\bf 323}, Springer: Berlin.

\bibitem{Oswald.P1999c}
Oswald P 1999,
{\em On the robustness of the {BPX}-preconditioner with respect to jumps 
in the coefficients}, Math. Comp., {\bf 68}, 633-650.

\bibitem{pazy} Pazy A 1983,
       {\em Semigroups of Linear Operators and Applications to Partial 
       Differential Equations}, New York: Springer.

\bibitem{reed} Reed M and Simon B, 1980, 1975, 1978, {\em
Methods of modern mathematical physics}, Volume I, II, IV,
Academic: New York.

\bibitem{weidmann} Weidmann J 1980, {\em Linear Operators in
 Hilbert Spaces}, Springer: New York.

\bibitem{XuZh:2008}
Xu J and Zhu Y 2008,
{\em Uniform convergent multigrid methods for elliptic problems with 
strongly discontinuous coefficients}, 
Math. Models Methods Appl. Sci., {\bf 18}, 77-105.

\bibitem{zhu:2008}
Zhu Y 2008,
{\em Domain decomposition preconditioners for elliptic equations with jump coefficients},
Numer. Linear Algebra Appl., {\bf 15}, 271-289.

\end{thebibliography}
\end{document}